\documentclass[a4paper,Haag duality]{mathscan}
\newtheorem{thm}{Theorem}[section] 
\newtheorem{pro}[thm]{Proposition}

\theoremstyle{definition}           
\newtheorem{rem}[thm]{Remark}       

\newcommand{\NI}{\noindent}

\newcommand{\bea}{\begin{eqnarray}}
\newcommand{\eea}{\end{eqnarray}}
\newcommand{\dsp}{\displaystyle}
\def \b #1 {\bf #1}
\newcommand{\IR}{\mathbb{R}}

\newcommand{\IM}{\mathbb{M}}
\newcommand{\IE}{\mathbb{E}}
\newcommand{\IC}{\mathbb{C}}
\newcommand{\ID}{\mathbb{D}}

\newcommand{\IZ}{\mathbb{Z}}

\newcommand{\cal}{\mathcal}

\newcommand{\cla}{{\cal A}}
\newcommand{\clm}{{\cal M}}

\newcommand{\clf}{{\cal F}}

\newcommand{\clh}{{\cal H}}
\newcommand{\clp}{{\cal P}}

\newcommand{\clb}{{\cal B}}

\newcommand{\clj}{{\cal J}}
\newcommand{\cln}{{\cal N}}
\newcommand{\cld}{{\cal D}}
\newcommand{\cll}{{\cal L}}
\newcommand{\clc}{{\cal C}}

\newcommand{\raro}{\rightarrow}

\newcommand{\vsp}{\vskip 1em}

\newcommand{\ul}{\underline}
\newcommand{\ol}{\overline}

\def \qed {\hfill \vrule height6pt width 6pt depth 0pt}
\newcommand{\be}{\begin{equation}}
\newcommand{\ee}{\end{equation}}
\newcommand{\ben}{\begin{eqnarray*}}
\newcommand{\een}{\end{eqnarray*}}
\begin{document}

\title{Translation invariant state and its mean entropy-II}

\author{ Anilesh Mohari }

\address{ The Institute of Mathematical Sciences, CIT Campus, Taramani, Chennai-600113 }

\email{anilesh@imsc.res.in}



\thanks{ }

\begin{abstract}

Let $\IM =\otimes_{n \in \IZ}\!M^{(n)}(\IC)$ be the two sided infinite tensor product $C^*$-algebra of $d$ dimensional matrices $\!M^{(n)}(\IC)=\!M_d(\IC)$ over the field of complex numbers $\IC$. Let $\omega$ be a translation invariant state of $\IM$. In a recent paper, we have proved that the mean entropy $s(\omega)$ is a complete invariant for certain classes of translation invariant state $\omega$ of $\IM$. In this paper, we have developed a general theory for dynamical entropy for an automorphism on an arbitrary $C^*$- or von-Neumann algebras based on repeated admissible measurement processes. In particular, we prove that dynamical entropy $h_{\omega}(\theta)$ for translation dynamics $(\IM,\theta,\omega)$ satisfies $s(\omega) \le h_{\omega}(\theta) \le 2s(\omega)$. In case $\omega$ is an infinite tensor product state of $\IM$ then $h_{\omega}(\theta)=s(\omega)$.    

\end{abstract}

\maketitle 

\section{ Introduction }

\bigskip
A measurement in classical dynamical systems gives a measurable partition
$\hat{\ul{\zeta}}=(\hat{\zeta}_i)$ of the configuration space ( probability measure space ) $(M,\clb,\mu)$.
For such a partition $\hat{\ul{\zeta}}=(\hat{\zeta}_i)$ of the configuration space, we can associate a family 
$\mu_{\hat{\zeta}_i}(E)={ \mu(E \cap \hat{\zeta}_i) \over \mu(E) }$ of probability measures and check that the Shannon information defined by 
\be 
H_{\mu}(\hat{\ul{\zeta}})=-\sum_i\mu(\hat{\zeta}_i) ln\mu(\hat{\zeta}_i)
\ee 
can be re-expressed as $\sum_i\mu(\hat{\zeta}_i)S(\mu_i,\mu)$, where $S(\mu_i,\mu)=\int {d\mu_i \over d\mu}ln({d\mu_i \over d\mu})d\mu$ is the Kullbeck-Liebler divergence or relative entropy between two probability measures $\mu_i$ and $\mu$. In other words, $H_{\mu}(\hat{\ul{\zeta}})$
can be interpreted as the average Kullbeck-Liebler divergence of the possible
final measures $(\mu_i)$ with respect to the initial probability measure $\mu$.
One of the important feature of the classical measurement is `the invariance properties' i.e. $\sum_i\mu(\hat{\zeta}_i)\mu_i=\mu$, which reflects that a measurement on and average does not disturb 
the classical system.  

\vsp 
Given two such partitions or measurements $\hat{\ul{\zeta}}=(\hat{\zeta}_i),\;\hat{\ul{\eta}}=(\hat{\eta}_j)$, we write 
$\hat{\ul{\zeta}} \circ \hat{\ul{\eta}} = (\hat{\zeta}_i \cap \hat{\eta}_j)$ as their joint measurements. One more feature in classical information theory is the following sub-additive property i.e. 
$$H_{\mu}(\hat{\ul{\zeta}} \circ \hat{\ul{\eta}} ) \le H_{\mu}(\hat{\ul{\zeta}}) + H_{\mu}(\hat{\ul{\eta}})$$
The above sub-additive property ensured existence for a notion for dynamical entropy $h_{\phi}(\theta,\ul{\hat{\zeta}})$ for an automorphism $\theta$ on $\clm=L^{\infty}(M,\clb,\mu)$ that preserves the measure $\mu$ and an invariance for the dynamics $h_{\mu}(\theta)$ defined by $$h_{\mu}(\theta)=\mbox{sup}_{\hat{\clp}}h_{\mu}(\theta,\ul{\hat{\zeta}}),$$
where sup is taken over all possible measurable partitions $\hat{\clp}$ of $M$. For further details we refer to original work 
described in details in [Si] and [CFS].

\vsp
In the quantum situation, we need to deal with a general von-Neumann sub-algebra $\clm$ of $\clb(\clh)$ that need not be isomorphic to the commutative algebra $L^{\infty}(M,\clb,\mu)$ and a normal state $\phi$ on $\clm$ replacing the role played by the probability measure $\mu$, where $\clb(\clh)$ is the algebra of bounded operator acting on a Hilbert space $\clh$. For an automorphism $\theta$ on $\clm$ preserving normal state $\phi$, we desires to introduce an invariance for $(\clm,\theta,\phi)$ as dynamical entropy $h_{\phi}(\theta)$ based on quantum mechanical measurement processes [Ne],[OP]. 

\vsp 
A quantum measurement, according to von-Neumann [Ne], give rises to a spectral resolution of identity operator i.e. a countable family of orthogonal projections $(\hat{\zeta}_i \in \clm)$ such that $\sum_i \hat{\zeta}_i = I$ and the final possible normal states are $\phi_i(x)= {1 \over \phi((\hat{\zeta}_i)}\phi(\hat{\zeta}_ix\hat{\zeta}_i)$ with a probability distribution given by $(\phi(\hat{\zeta}_i))$. The major difficulties with this family of measurements or partitions of unity are the following:

\NI (a) for two such given partitions $\hat{\ul{\zeta}},\hat{\ul{\eta}}$ of unity, there is no natural 
meaning for joint measurements or partition of unity unless two partitions mutually commutes 
$\hat{\zeta}_i \hat{\eta}_j=\hat{\eta}_j \hat{\zeta}_i$. 

\NI (b) In a quantum situation, there are few projection valued measurements that keep the state $\phi$ invariant i.e. the following equality may not hold: 
$\sum_i\phi(\hat{\zeta}_i)\phi_i(x) = \phi(x) \; \forall x \in \clm$.

\vsp 
These motivate us to work beyond the framework of projection valued von-Neumann measurements. By accepting the theme that a measurement is an irreversible process in quantum situation, we consider a more general measurement process 
given by a family of completely positive maps $\ul{\zeta}=(\zeta_i)$ on $\clm$ such that $\zeta=\sum_i \zeta_i$ is a unital map on $\clm$. Such a family $\ul{\zeta}=(\zeta_i)$ of completely positive maps will be referred as {\it measurement} 
or {\it a partition of unity} in $\clm$. For any two partitions $\ul{\zeta}$ and $\ul{\eta}$ of unity in $\clm$, we set partition $\ul{\zeta} \circ \ul{\eta}$ of unity in $\clm$ by $\ul{\zeta} \circ \ul{\eta} = (\zeta_i\eta_j)$. We say 
$\ul{\zeta}$ commutes with $\ul{\eta}$ if $\zeta_i \eta_j = \eta_j \zeta_i$. A partition of unity $\ul{\zeta}=(\zeta_i)$ in $\clm$ is called {\it projection valued measurement} if $\zeta_i(x)=\hat{\zeta}_ix\hat{\zeta}_i$ for all $x \in \clm$ for a family $\hat{\ul{\zeta}}=(\hat{\zeta}_i)$ of orthogonal projections in $\clm$.  

\vsp 
For (1), we adopt the quantum mechanical Kullbeck-Liebler ( Araki's  relative entropy ) divergence of the average final state with the final possible state. In other-words we define
\be 
H_{\phi}(\ul{\zeta})=\sum_i\phi(\zeta_i(I))S({ \phi \circ \zeta_i \over \phi(\zeta_i(I)) }, 
\phi \circ \zeta) 
\ee 

\vsp 
For any two partitions $\ul{\zeta}$ and $\ul{\eta}$ of unity in $\clm$, we prove that
\be
H_{\phi}(\ul{\zeta} | \ul{\eta} \circ \ul{\beta}) 
\le H_{\phi}(\ul{\zeta}|\ul{\eta})
\ee
where 
\be 
H_{\phi}(\ul{\zeta}|\ul{\eta})=
H_{\phi}(\ul{\zeta} \circ \ul{\eta}) - H_{\phi \circ \zeta}(\ul{\eta})
\ee

\vsp 
The basic inequality (3) ensures an existence theorem of a dynamical entropy $h_{\phi}(\theta,
\ul{\zeta})$ defined by
\be 
h_{\phi}(\theta,\ul{\zeta})= \mbox{limit}_{n \raro \infty}H_{\phi}(\ul{\zeta}|\ul{\zeta}^{-}_{(n)})
\ee 
for any partition $\ul{\zeta}=(\zeta_i)$ of unity and a $\phi$ preserving automorphism
$\theta$, where $\ul{\zeta}^{-}_{(n)}=\theta^{-1}(\ul{\zeta})\circ \theta^{-2}(\ul{\zeta})\circ ....
\circ \theta^{-n}(\ul{\zeta})$
and $\theta(\ul{\zeta})=\theta\ul{\zeta}\theta^{-1}.$ It is also simple to check that 
$h_{\phi}(\theta,\ul{\zeta})=h_{\phi}(I,\theta \ul{\zeta})$ and  
$\theta \circ \ul{\zeta} = (\theta\zeta_i)$ is also a partition of unity in $\clm$ 
once $(\zeta_i)$ is so. Thus the class of measurements is too large to destroy the characteristic for 
a particular dynamics $(\theta)$. In mathematical term, we note that 
$\mbox{sup}_{\ul{\zeta} \in \clp}h_{\phi}(\theta,\ul{\zeta})$ is independent of $\theta$, if sup is taken 
over all possible partitions $\clp$ of unity in $\clm$. A partition $\ul{\zeta}=(\zeta_i)$ is called $\phi$-{\it invariant} 
for the state $\phi$ if $\phi \circ \zeta = \phi$. Note also that $\mbox{sup}_{\ul{\zeta} \in \clp_{\phi}}h_{\phi}(\theta,\ul{\zeta})$ 
is as well independent of $\theta$, if sup is taken over all possible $\phi$-invariant partitions $\clp_{\phi}$ of unity in $\clm$.

\vsp
Thus we need to look at a smaller class of measurements. We say a partition or measurement $(\zeta_i)$ in $\clm$ is of {\it zero mean entropy for } $(\clm,I,\phi)$ if $h_{\phi}(I,\ul{\zeta})=0$. For a partition $(\hat{\zeta}_i)$ of identity into orthogonal projections on $\clh$, we define partition given by the family of completely positive maps $\zeta_i$ by $\zeta_i(x)=\hat{\zeta}_ix\hat{\zeta}_i,\;\forall x \in \clm$. It is obvious that $h_{\phi}(I,\ul{\zeta})=0$ for such a partition $\ul{\zeta}=(\zeta_i)$ of unity. However, such a von-Neumann measurement in general need not be invariant for a normal state $\phi$. Nevertheless the class of invariant admissible measurements are not small. As an example, we may recall [Mo5] translation automorphism $\theta$ on the two-sided quantum spin chain $\IM=\otimes_{k \in \IZ}\!M^{(k)}_d(\IC)$ with an invariant state $\omega$. Proposition 5.4 in [Mo5] says that there exists an automorphism $\alpha$ commuting with $\theta$ such that $\omega \alpha \IE_0= \omega \alpha$, where $\IE_0$ is the trace preserving norm one projection on the maximal abelian $C^*$-sub-algebra $\ID_e=\otimes_{k \in \IZ} \ID^{(k)}_e(\IC)$ with $\ID^{(k)}_e(\IC)$ be the diagonal matrices with respect to an 
orthonormal basis $e=(e_j)$ of $\IC^d$. A natural partition of unity in $\clm$ can be described by $\eta_j(x)= \pi_{\omega}(|e_j \rangle \langle e_j|) \IE_0(x) \pi_{\omega}(|e_j \rangle \langle e_j|)$ for $x \in \clm$. We also note that $(\alpha \eta_j)$ is a $\omega$-invariant admissible partition of unity in $\clm$ 

\vsp 
A partition $\ul{\zeta}$ of unity in $\clm$ is called {\it von-Neumann for the automorphism} $\theta$ on $\clm$ if the family $\{(\theta^k(\zeta_i)): k \in \IZ \}$ of partitions are mutually commuting. One important difference here, we are not demanding the family $(\zeta_i)$ itself to be mutually commuting. We finally say a measurement $\ul{\zeta}=(\zeta_i)$ is {\it admissible for} $(\clm,\theta,\phi)$ if $\ul{\zeta}$ is a von-Neumann measurement for $(\clm,\theta^k,\phi)$ for some $k \ge 1$. In such a case we define 
\be 
h_{\phi}(\theta,\ul{\zeta})= \mbox{sup}_{ k } {1 \over k} h_{\phi}(\theta^k,\ul{\zeta})
\ee
where sup is taken over all possible values for $k$ for which $\ul{\zeta}$ is a von-Neumann partitions 
in $(\clm,\theta^k,\phi)$. We also note that partition $(\eta_j)$ of $\pi_{\omega}(\IM)''$ described above is admissible for the translation dynamics $\theta$ on $\pi_{\omega}(\IM)''$.  

\vsp
We define {\it quantum dynamical entropy} $h_{\phi}(\theta)$ by taking sup of $h_{\phi}(\theta,
\ul{\zeta})$ over all admissible $\phi$-invariant measurements $\clp_{\phi,a}$ i.e. 
\be 
h_{\phi}(\theta)=\mbox{sup}_{\ul{\zeta} \in \clp_{\phi,a} }h_{\phi}(\theta,\ul{\zeta})
\ee
Thus $h_{\phi}(\theta)$ is an invariance for the dynamics. A question that is central now: given an automorphism $\theta$ on a von-Neumann algebra $\clm$ with a faithful normal invariant state $\phi$, 
is there enough choices for $\phi-$invariant $\theta$-admissible partition of unity in $\clm$ to characterise the dynamics $(\clm,\theta,\phi)$ for a given value of $h_{\phi}(\theta)$? This problem in the classical framework gave remarkable results [Or1] once restricted to certain class of automorphisms [Or2].  

\vsp 
In the last section, we consider the translation dynamics $(\IM,\theta,\omega)$ and proved that $s(\omega) \le h_{\omega}(\theta) \le 2 s(\omega)$. In case $\omega$ is an infinite tensor product state of $\IM$ then $h_{\omega}(\theta)=s(\omega)$. 
As of now, we do not have a counter example to suggest $h_{\omega}(\theta) \neq s(\omega)$. We postpone dealing with other examples of translation invariant states that arises in non-commutative quantum dynamics to a forth coming paper [Mo2].   

\vsp
\section{ Araki's relative entropy and quantum information: }

\vsp 
Let $\clm$ be a von-Neumann algebra and $\omega$ be 
a faithful normal state. Without loss
of generality let also $(\clm,\omega)$ be in the standard form
$(\clm,J_{\omega},{\cal P}_{\omega},\zeta_{\omega})$
[BrR] where $\zeta_{\omega}
\in \clh$, a cyclic and separating vector for $\clm$, so
that $\omega(x)= <\zeta_{\omega},x\zeta_{\omega}>$ and the closer of the
closable operator $S^0_{\omega}:x\zeta_{\omega}
\raro x^*\zeta_{\omega}, S_{\omega}$ possesses a polar decomposition
$S_{\omega}=J_{\omega}\Delta^{1/2}_{\omega}$ with the self-dual positive
cone $\clp_{\omega}$ as the closure of $\{ J_{\omega}xJ_{\omega}x\zeta_{\omega}:x \in \clm \}$ in
$\clh$. Tomita's [BrR] theorem says that
$\Delta_{\omega}^{it}\clm\Delta_{\omega}^{-it}=\clm,\;t
\in \IR$ and $J_{\omega} \clm J_{\omega}=\clm'$, where $\clm'$ is the
commutant of $\clm$. We define the modular automorphism group
$\sigma^{\omega}=(\sigma^{\omega}_t,\;t \in \IR )$ on $\clm$
by
$$\sigma^{\omega}_t(x)=\Delta_{\omega}^{it}x\Delta_{\omega}^{-it}.$$
Furthermore for any normal positive functional $\psi$ on $\clm$
there exists a unique vector $\zeta_{\psi} \in {\cal P}_{\omega}$ so that $\psi(x)=
<\zeta_{\psi},x\zeta_{\psi}>$.

\smallskip
Following Araki [Ar,OP] we define the relative entropy $S(\psi_1,\psi_2)$
for two normal positive functionals $\psi_1,\psi_2$ on $\clm$ where
$\psi_1(x)=<\zeta_{\psi_1},x \zeta_{\psi_1}>$ and $\psi_2(x)=
<\zeta_{\psi_2},x\zeta_{\psi_2}>$ for some unique 
$\zeta_{\psi_1},\zeta_{\psi_2} \in {\cal P}$. The closer of the
closable operator 
$$S^0_{\psi_2,\psi_1}:x\zeta_{\psi_2}+z
\raro s^{\clm}(\zeta)x^*\zeta_{\psi_1},\;x \in \clm, s^{\clm'}z=0$$
defined on $\clm\zeta_{\psi_2} + (I-s^{\clm'}(\zeta_{\psi_2}))\clh$
where $s^{\clm}(\zeta_{\psi_2})$ is the projection from $\clh$ to $\{\ol{\clm'\zeta_{\psi_2}
\}}$
, $S(\zeta_{\psi_2},\zeta_{\psi_1})$ possesses a polar decomposition
$S_{\psi_2,\psi_1}=J_{\psi_2,\psi_1}\Delta^{1/2}_{\psi_2,\psi_1} $. So by
definition, $\Delta_{\omega,\omega}=\Delta_{\omega}$. The Araki's relative entropy
is defined by
$$S(\psi_2,\psi_1)=\{ \begin{array}{lll}-\langle \zeta_{\psi_2},ln(\Delta_{\psi_2,\psi_1})
\zeta_{\zeta_2}\rangle,\;
&if\; \psi_2 << \psi_1 \\ \infty,\;&\;otherwise \end{array}.$$
where $\psi_2 << \psi_1$ means that $\psi_2(x^*x)=0$ implies $\psi_1(x^*x)=0$ for
$x \in \clm$. We recall in the following proposition few well-known
properties of relative entropy.

\vsp
\begin{pro} The relative entropy of two positive functional satisfies the following relations;

\NI (a) $S(\mu\omega,\lambda \phi) = \mu S(\omega,\phi)-\mu\omega(I)(log\lambda -log \mu)$ for any 
$\mu,\lambda \ge 0$;

\NI (b) $S(\omega,\phi) \ge \omega(I)(log \omega(I) - log\phi(I))$

\NI (c) Jointly convex: $S(\psi,\phi)$ is jointly convex, i.e.
$S(\lambda \psi_1 + \mu \psi_2, \lambda \phi_1 + \mu \phi_2) \le
\lambda S(\psi_1, \phi_1) + \mu S(\psi_2, \phi_2)$, where $\mu +\lambda=1,\mu
,\lambda \ge 0$.

\NI (d) Ullhamm's monotonicity: For any Schwartz type positive unital normal
map $\tau$ from $\clm$ to $\cln$, $$S(\psi_2\tau,\psi_1\tau) \le 
S(\psi_2,\psi_1).$$

\NI (e) Lower semi-continuity: If $\mbox{lim}_{n \raro \infty}\psi_n=\psi$
and $\mbox{lim}_{n \raro \infty}\phi_n=\phi$ in weak$^*$ topology, then 
$S(\psi,\phi) \le \mbox{lim inf}_{n \raro \infty} S(\psi_n,\phi_n)$. 
Moreover; if there exists a positive number $\lambda$ satisfying $\psi_n 
\le \lambda \phi_n$, then $\mbox{lim}_{n \raro \infty}
S(\psi_n,\phi_n)=S(\psi,\phi)$.

\NI (f) Donald's identity: For any family of normal positive functional
$(\psi_i,\;i \ge 1)$ and $\phi$ we have
$$S(\psi,\phi) + \sum_i S(\psi_i,\omega) = \sum_i S(\psi_i,\phi)$$
where $\sum_i\psi_i=\psi$.

\NI (g) For two normal state $\psi_1,\psi_2$, we have the the lower lower bound: 
$||\psi_1-\psi_2||^2/2 \le S(\psi_2,\psi_1)$.

\NI (h) Let $\clm_0$ be a von-Neumann sub-algebra of $\clm$. Assume that there exists a faithful normal norm one projection $\IE_0$ from $\clm$ onto $\clm_0$. If $\psi_0$ and $\phi$ are normal states of $\clm_0$ and $\clm$ respectively, then 
\be 
S(\phi, \psi_0 \circ \IE_0)= S(\phi|\clm_0,\psi_0)+ S(\phi,\phi \circ \IE_0)
\ee    
\end{pro} 

\vsp
\begin{proof}
We refer to the monograph [OP]. 
\end{proof} 

\bigskip
Let $(\Omega,\clf,\mu)$ be a probability space. A classical dynamics $(\Omega,\theta,\mu)$ is a triplet, 
where $\theta:\Omega \raro \Omega$ is bi-measurable one to one and onto map modulo a $\mu$-null set and 
$\mu \circ \theta^{-1} = \mu$ on measurable subsets of $\Omega$. This classical dynamics is equivalently described by $(\clm_0,\theta,\phi_{\mu})$, where $\clm_0=L^{\infty}(\Omega,\clf,d\mu)$ and $\theta:\clm_0 \raro \clm_0$ be the automorphism defined by $\theta(f)= f \circ \theta$ for all $f \in \clm_0$ with normal invariant state $\phi_{\mu}(f)=\int f d\mu$. A measurement in classical dynamics $(\Omega,\theta,\mu)$ gives a measurable partition $\ul{\zeta}=(\zeta_i)$ of the measure space $(\Omega,\clb,\mu)$. For such a measurable partition $\ul{\zeta}=(\zeta_i)$, we can associate a family of probability measures $\mu_{\zeta_i}(E)={\mu(E \cap \zeta_i) \over \mu(E)}$ provided $\mu(\zeta_i) > 0$ and check that Shannon information [Pa] defined by 
$$H_{\mu}(\ul{\zeta})=-\sum_i\mu(\zeta_i) ln\mu(\zeta_i)$$ 
can be re-expressed as 
$$H_{\mu}(\ul{\zeta})=\sum_i\mu(\zeta_i)S(\mu_i,\mu)$$ 
where $S(.,.)$ is the Kullbeck-Liebler relative entropy. In other words $H_{\mu}(\ul{\zeta})$
can be interpreted as the average Kullbeck-Liebler divergence of the final measures $(\mu_i)$ with respect to the initial measure $\mu$. Given two such measurable partitions or measurements $\ul{\zeta}=(\zeta_i),\;\ul{\eta}=(\eta_j)$, we write $\ul{\zeta} \circ \ul{\eta} = (\zeta_i \cap \eta_j)$ as their joint measurements. Thus two such measurements commute and the class of measurements admits invariance property i.e. $$\sum_i\mu(\zeta_i)\mu_i=\mu,$$ 
which reflects the fact that a classical measurement does not disturb the system. Furthermore, the classical information admits the following sub-additive property, i.e. 
$$H_{\mu}(\ul{\zeta} \circ \ul{\eta} ) \le H_{\mu}(\ul{\zeta}) + H_{\mu}(\ul{\eta})$$

\vsp
A quantum dynamics is a triplet $(\clm,\theta,\phi)$, where $\clm$ is a von-Neumann algebra acting on a complex separable Hilbert space and $\theta$ is a $*$ automorphism with a normal invariant state $\phi$ of $\clm$. 

\vsp 
We fix any two von-Neumann algebras $\clm$ and $\cln$. A family of completely positive maps $\underline{\zeta}=(\zeta_i:\clm \raro \cln)$ is called {\it partition} of unity for $\cln$ if $\sum_i\zeta_i(I_{\clm})=I_{\cln}$, where $I_{\clm}$ and $I_{\cln}$ are unit elements in $\clm$ and $\cln$ respectively. We also set unital completely positive map $\zeta=\sum_i\zeta_i:\clm \raro \cln$. 
So for any positive normal state $\phi$ on $\cln$ as an input or initial state,
$\phi_i=\phi(\zeta_i(I_{\clm}))^{-1}\phi \circ \zeta_i$ is the final or output state
on $\clm$ with probability $\phi(\zeta_i(I_{\clm}))$. We define {\it quantum information} associated with a 
measurement $\ul{\zeta}$ with input state $\phi$ on $\cln$ and out put state $\phi \zeta$ 
by 
\be 
H_{\phi}(\ul{\zeta})=\sum_i\phi(\zeta_i(I_{\clm}))S({ \phi \circ \zeta_i \over \phi(\zeta_i(I)) }, 
\phi \circ \zeta) 
\ee 
where $(\zeta_i)$ is a family of completely positive maps from $\clm$ to $\cln$ 
so that $\zeta=\sum_i\zeta_i$ is a unital map from $\clm$ to $\cln$. 

\vsp
We consider the von-Neumann algebra $\clm_{\zeta} = \oplus_{i \in \zeta}
\clm_i$ where each $\clm_i$ are copies of the von-Neumann $\clm$
and two states $\phi^1_{\zeta}(a)=\dsp \sum_{i \in \zeta}
\phi(\zeta_i(a_i)) $
and $\phi^2_{\zeta}(a)=
\sum_{i \in \zeta}\phi(\zeta_i(I_{\clm}))\phi(\zeta(a_i))$.
By rescaling and additivity property of relative entropy we verify
that
\be
H_{\phi}(\ul{\zeta})=S(\phi^1,\phi^2).
\ee
and also
\be
H_{\phi}(\ul{\zeta})=H^c_{\phi}(\ul{\zeta}) + H^q_{\phi}(\ul{\zeta})
\ee
where 
$H^c_{\phi}(\ul{\zeta})= \sum_{i \in \zeta}-\phi(\zeta_i(I_{\clm}))ln(\phi(\zeta_i(I_{\clm}))$ 
and 
$H^q_{\phi}(\ul{\zeta})= \sum_{i \in \zeta} S(\phi \circ \zeta_i, \phi \circ \zeta)$

\vsp
For any two partitions $\underline{\zeta} =(\zeta_i:\clm \raro \cln)$ and $\underline{\eta}=(\eta_j:\clc \raro \clm)$, we set partition $\underline{\zeta \circ \eta}
= ( \zeta_i\eta_j:\clc  \raro \cln \;\;i \in \zeta,j \in \eta )$. 

\vsp
\begin{pro} 
For any partition $\ul{\zeta}=(\zeta_i)$ the map $\phi \raro H_{\psi}(\ul{\zeta})$ is 
continuous. Also the map $\ul{\zeta} \raro H_{\psi}(\ul{\zeta})$ is continuous in weak$^*$ topology. Furthermore for any two partitions $\underline{\zeta}:\clm \raro \cln,\;
\underline{\eta}: \clc \raro \clm$ of unities in $\cln$ and $\clm$ respectively we have 
\be 
H_{\phi \circ \zeta}(\ul{\eta}) \le H_{\phi}(\ul{\zeta} \circ \ul{\eta}) 
\ee
and 
\be
H_{\phi}(\ul{\zeta} \circ \ul{\eta}) \le H_{\phi \circ \zeta}(\ul{\eta}) + H_{\phi}(\ul{\zeta})
\ee
Furthermore, the map $\phi \raro H^q_{\phi}(\ul{\zeta} \circ \ul{\eta})-
H^q_{\phi \circ \zeta}(\ul{\eta})$ is convex.
\end{pro}

\vsp
\begin{proof} 
First part is an easy consequence of joint continuity property (e) in Proposition 2.1. We claim the 
following two inequalities:
$$
-\sum_{i,j}\phi(\zeta_i\eta_j(I))ln(\phi(\zeta_i\eta_j(I_{\clc})) \le -\sum_{i \in \zeta}\phi(\zeta_i(I_{\clm}))ln(\phi(\zeta_i(I_{\clm}))) 
$$
\be 
- \sum_{j \in \eta}\phi(\zeta\eta_j(I_{\clc}))ln(\phi(\zeta\eta_j(I_{\clc})))
\ee
and
\be
\sum_{i \in \zeta, j \in \eta} S(\phi\zeta_i\eta_j,\phi \zeta\eta) \le \sum_{j \in \eta}S(\phi
\zeta\eta_j,\phi\zeta\eta) + \sum_{i \in \zeta}S(\phi \zeta_i,\phi \zeta)
\ee
For the first inequality (14), we note that both $\{ \phi(\zeta_i\eta_j(I_{\clc})):\;i \in \zeta,\;j \in \eta \}$ and $\{\phi(\zeta_i(I_{\clm})\phi(\zeta\eta_j(I_{\clc})):i \in \zeta, \; j \in \eta \} $ are probability
measure on the index set $\zeta \times \eta$. The inequality is nothing but the well known statement that relative entropy of any two probability measure is always non-negative. For the second
inequality we first appeal to Donald's identity to check that
\be
\sum_{i \in \zeta, j \in \eta} S(\phi\zeta_i\eta_j,\phi \zeta\eta)= 
\sum_{j \in \eta}S(\phi
\zeta\eta_j,\phi\zeta\eta) + \sum_{i \in \zeta, j \in \eta}S(\phi \zeta_i \eta_j,\phi \zeta 
\eta_j)
\ee
Now for each $i \in \zeta$ we check that $\sum_{j \in \eta}S(\phi \zeta_i\eta_j,\phi \zeta\eta_j)
=S(\phi\zeta_i\IE_{\eta},\phi\zeta\IE_{\eta})$ where $\IE_{\eta}:\clc_{\eta} \raro \clm$ is
the unital Schwartz type map defined by $$\IE_{\eta}((c_j))= \sum_{j \in \eta}\eta_j(c_j).$$
Thus by Ullhamm's monotonicity we conclude the second inequality (15).
By combing (14) and (15), we conclude the proof for (13). The convex property of the map $\phi \raro H^q_{\phi}(\ul{\zeta}|\ul{\eta})=H^q_{\phi}(\ul{\zeta}\circ \ul{\eta})-H^q_{\phi}(\ul{\eta})$ 
follows by the convexity property of relative entropy once we appeal to the identity 
(16). 
\end{proof}

\vsp
The following proposition is an useful generalization of Proposition 3.1.
\vsp
\begin{pro} 
For any three partitions $\underline{\zeta}:\cln
\raro \clm,\;\underline{\eta}:\clc \raro \cln,\;\underline{\beta}:\cld \raro \clc$ of unities:
\be
0 \le H_{\phi}(\zeta \circ \eta \circ \beta)-H_{\phi \circ \zeta}(\eta \circ \beta) 
\le H_{\phi}(\zeta \circ \eta) -H_{\phi \circ \zeta}(\eta)
\ee
for any normal state $\phi$ on $\clm$.
\end{pro}

\vsp
\begin{proof} 
Essential steps are same as in Proposition 3.1. By Donald's
identity and Ullhamm's monotonicity we check 
that 
$$\sum_{i \in \zeta,j \in \eta, k \in \beta}S(\phi \zeta_i\eta_j\beta_k,\phi
\zeta\eta\beta)-\sum_{j \in \eta, k \in \beta}S(\phi\zeta\eta_j\beta_k,
\phi\zeta\eta\beta)  
\le \sum_{i \in \zeta,j \in \eta}S(\phi\zeta_i\eta_j,\phi\zeta\eta) $$
$$ - \sum_{j \in \eta} S(\phi\zeta\eta_j,\phi\zeta \eta_j) $$
Also by Jensen's inequality we check that 
$$\sum_{i \in \zeta,j \in \eta, k \in \beta}\phi \zeta_i\eta_j\beta_k(I_{\cld})ln
{\phi\zeta_i\eta_j\beta_k(I_{\cld}) \over \phi \zeta\eta_j\beta_k(I_{\cld}) } 
\le \sum_{i \in \zeta,j \in \eta}\phi\zeta_i\eta_j(I_{\clc})ln{\phi\zeta_i\eta_j(I_{\clc}) 
\over \phi\zeta \eta_j(I_{\clc})} $$ 
\end{proof}

\vsp 
In case $\clm=\cln$, such a family $\ul{\zeta}=(\zeta_i)$ will be referred as a {\it measurement} in $\clm$. For any two measurements $\zeta,\eta$ in $\clm$, we set measurement $\zeta \circ \eta$ in $\clm$ by $\zeta \circ \eta =(\zeta_i \circ \eta_j)$. 
Two measurements $\ul{\zeta}$ and $\ul{\eta}$ are called {\it mutually commuting} if $\ul{\zeta} \circ \ul{\eta}=\ul{\eta} \circ \ul{\zeta}$, where equality as sets. In particular, in such a case we have $H_{\phi}(\ul{\zeta} \circ \eta)=H_{\ul{\eta} \circ \ul{\zeta}}$. We also say a measurement $(\zeta_i)$ is {\it invariant} for the state $\phi$ if $\phi \circ \zeta = \phi$. Thus $\ul{\zeta} \circ \ul{\eta}$ is also $\phi$-invariant if $\ul{\zeta}$ and $\ul{\eta}$ are so. 

\vsp  
We have the following important observation:

\vsp 
\begin{pro} 
Let $\ul{\zeta},\ul{\eta}$ and $\ul{\beta}$ be three $\phi$-invariant partitions of unity in $\clm$ such that either $\ul{\zeta} \circ \ul{\eta} = \ul{\eta} \circ \ul{\zeta}$ or $\ul{\eta} \circ \ul{\beta}= \ul{\beta} \circ \ul{\eta}$. Then 
\be 
H_{\phi}(\ul{\zeta} \circ \ul{\eta} \circ \ul{\beta}) - H_{\phi}(\ul{\eta} \circ \ul{\beta})
\le H_{\phi}(\ul{\zeta} \circ \ul{\beta}) - H_{\phi}(\ul{\beta})
\ee 
\end{pro}

\vsp 
\begin{proof} 
Both the situations are simple consequences of Proposition 2.3. 
\end{proof} 

\vsp 
\section{Dynamical Entropy for $C^*$ or $W^*$-systems:} 

\bigskip
Now we aim to develop, a quantum mechanical analogue of dynamical entropy introduced by Kolmogorov and Sinai [See e.g. Pa] based on von-Neumann measurement process. We prove first a general existence theorem and discuss its analytical properties. We fix a von-Neumann algebra $\clm$ and assume that it is in standard form. We adopt the same notation of section 2 for measurements and channels assuming that both input and output algebras $\clm$ and $\cln$ are same.

\vsp 
Let $\theta$ be a $*$-automorphism on $\clm$ and $\phi$ is a normal state, invariant for $\theta$. For a partition $\ul{\zeta}=(\zeta_i)$ of unity in $\clm$ and an automorphism $\theta:\clm \raro \clm$, we set partition of unity 
$\theta(\ul{\zeta})=(\theta \zeta_i \theta^{-1})$ in $\clm$. Also note that $\ul{\theta} \circ \ul{\zeta}=(\theta \zeta_i)$ 
is also a partition of unity in $\clm$. The partitions $\theta(\ul{\zeta})$ and $\ul{\theta} \circ \ul{\zeta}$ are $\phi$-invariant if $\ul{\zeta}$ and $\theta$ are so. For each $n \ge 1$, we set partition of unity $\ul{\zeta}_{(n)}= \theta^{n-1}(\ul{\zeta}) \circ \theta^{n-2}(\ul{\zeta}) \circ.....\circ \ul{\zeta}$ in $\clm$ and $\ul{\zeta}^{-}_{(n)}= \theta^{-1}(\ul{\zeta}) \circ \theta^{-2}(\ul{\zeta})... \circ \theta^{-n+1}(\ul{\zeta}) \circ \theta^{-n}(\ul{\zeta})$. So $\zeta_{(n)}=\theta^{n-1}(\zeta) \circ ...\theta(\zeta) \circ \zeta$ is an unital map on $\clm$.

\vsp
\begin{pro} 
For a unital $*$-automorphism $\theta$ on $\clm$ with
$\phi=\phi \circ \theta$   
\be
H_{\phi}(\ul{\zeta})=H_{\phi}(\theta(\ul{\zeta}))
\ee
for any partition $\zeta$, where $\theta(\ul{\zeta})_i=\theta \circ \zeta_i \circ \theta^{-1}.$
\end{pro}

\vsp
\begin{proof}   
Since the relative entropy of any two normal states remain invariant i.e. 
$S(\psi_1,\psi_2)=S(\theta\psi_1\theta^{-1},\theta\psi_2\theta^{-1})$, 
by an automorphism, the equality is immediate. 
\end{proof} 

\vsp
\begin{pro} 
For a countable partitions $\underline{\zeta}$ with $H_{\phi}(\underline{\zeta}) < \infty$,
\be 
h_{\phi}(\theta,\underline{\zeta})=\mbox{lim}_{n \raro \infty}H_{\phi}
(\underline{\zeta}|\underline{\zeta}^-_{(n)})
\ee
exists, where $\underline{\zeta}^-_{(n)}= \theta^{-1}(\underline{\zeta})\circ \theta^{-2}(\underline{\zeta}) \circ ..
\circ \theta^{-n}(\underline{\zeta})$. The map $(\phi,\ul{\zeta}) \raro h_{\phi}(\theta,\ul{\zeta})$ is upper semi-continuous 
in each variables in the Bounded Weak topology.
\end{pro}

\vsp 
\begin{proof}
We set $a_{\phi,n}(\ul{\zeta}) = H_{\phi}(\theta^n(\underline{\zeta})|\theta^{n-1}(\underline{\zeta})\circ
..\circ \underline{\zeta})=H_{\phi}(\ul{\zeta}|\ul{\zeta}^{-}_{(n)})$. By Proposition 3.1 and Proposition 4.1 we 
check that $ 0 \le a_{\phi,n+1}(\ul{\zeta}) \le a_{\phi,n}(\ul{\zeta}) \le
H_{\phi}(\underline{\zeta})$, 
thus the $\mbox{limit}_{n \raro \infty}a_{\phi,n}$ exists. By Proposition 2.1 (e), the map $(\phi,\ul{\zeta}) \raro a_{\phi,n}(\ul{\zeta})$ is continuous in Arveson's Bounded Weak topology [Pau]. Thus the last statement follows as the liming function is the inf over continuous functions. 
\end{proof} 

\vsp
Since $a_n(\ul{\zeta})$ are monotonically decreasing, their limit $h_{\phi}(\theta,\ul{\zeta})$
is also same as $\mbox{limit}_{n \raro \infty}{1 \over n} \sum_{1 \le k \le n}
a_k(\ul{\zeta}) $. So if we interpret $a_n(\zeta)$ as information gain in the 
nth step, $h_{\phi}(\theta, \ul{\zeta})$ is indeed the average information gained 
for large $n$. It is simple to check that for a partition $\ul{\zeta}$, 
$$h_{\phi}(\theta,\ul{\zeta}) = h_{\phi}(I,\ul{\theta} \circ \ul{\zeta}),$$
where $I$ is the identity map on the von-Neumann algebra $\clm$ and $\ul{\theta}=(\theta)$. 
Thus in case we intend to define dynamical entropy $h_{\phi}(\theta)$ of 
$(\clm,\theta,\phi)$ to be the supremum over all the partitions, i.e. 
$\{\ul{\zeta}=(\zeta_i:1 \le i \le d)$ then the value $h_{\phi}(I)$ will 
be independent of $\theta$. So $h_{\phi}(I)$ is an invariance for the von-Neumann 
algebra $\clm$ possibly with values infinite. 

\vsp 
In case $\phi \circ \zeta=\phi$, 
\be 
h_{\phi}(\theta,\ul{\zeta})=\mbox{limit}_{n \raro \infty} { 1 \over n}H_{\phi}(\theta^{n-1}(\ul{\zeta})\circ ...\ul{\zeta})
\ee 
For any $\phi$-invariant partition of unity in $\clm$, the partition $\theta \circ \ul{\zeta}$ is also 
a $\phi$-invariant and so $h_{\phi}(I)$ will be independent if sup is taken over all possible $\phi$-invariant partition and these non-negative numbers are invariance for the von-Neumann algebra $\clm$. It also suggest that we need to reduce the class of measurements if we want to extract information about a given automorphism $\theta$ on $\clm$. 

\vsp 
The major difficulties that we face while dealing with the class of partitions of unity in the non-commutative framework are the following:

\vsp 
\NI (a) For two given partitions $\ul{\zeta},\ul{\eta}$ of unity in $\clm$ arises from two measurements, 
there is no natural meaning to joint measurements unless we have $\ul{\zeta} \circ \ul{\eta} = \ul{\eta}
\circ \ul{\zeta}$ as sets. 

\vsp 
\NI (b) In quantum situation, a measurement $\ul{\zeta}$ need not keep the state $\phi$ invariant in general though desirable in the theory of repeated measurement proposed by von-Neumann [Ne].  

\vsp 
\begin{pro} 
Let $\ul{\zeta}=(\zeta_i)$ and $\ul{\eta}=(\eta_j)$ be two measurements in $\clm$. Then 
\be 
H_{\phi}(\ul{\zeta} \circ \ul{\eta}) \ge H_{\phi}(\ul{\zeta} \circ \eta)
\ee
\end{pro} 

\vsp 
\begin{proof} 
By the join convexity property of relative entropy: 
$$H_{\phi}(\ul{\zeta} \circ \ul{\eta}) = \sum_{i,j} \phi(\zeta_i\eta_j(I))  S( { 1 \over \phi(\zeta_i\eta_j(I))} \phi \zeta_i \eta_j, \phi \zeta \eta)$$
$$ \ge \sum_i \phi(\zeta_i(I)) S( {1 \over \phi \zeta_i(I) } \phi \zeta_i \eta, \phi \zeta \eta)$$ 
$$ = H_{\phi}(\ul{\zeta} \circ \eta)$$
\end{proof} 

\vsp 
A partition $\ul{\zeta}$ of unity in $\clm$ is called {\it von-Neumann measurement } 
for $(\clm,\theta,\phi)$ if the family $\{\theta^k(\ul{\zeta}):k \in \IZ \}$ of partitions of unity 
in $\clm$ are mutually commuting. In particular, $\theta^n(\ul{\zeta})$ and $\ul{\zeta}_{(n)}$ are 
mutually commuting for each $n \ge 1$. Note that we are not demanding the family $(\zeta_i)$ itself to be commutative.

\vsp 
For any two $\phi$-invariant partitions $\ul{\zeta}$ and $\ul{\eta}$ of unity in $\clm$ and an automorphism $\theta$ on $\clm$, we set notation 
$$H_{\phi}(\ul{\zeta}_{(n)}| \ul{\eta}_{(n)})= H_{\phi}( \theta^{n-1}(\ul{\zeta}) \circ \theta^{n-2}(\ul{\zeta}) \circ... \circ \theta(\ul{\zeta}) \circ \ul{\zeta} \circ 
\theta^{n-1}(\ul{\eta}) \circ \theta^{n-2}(\ul{\eta}) \circ ... \circ \theta(\ul{\eta}) \circ \ul{\eta})$$
$$ 
- H_{\phi}(\theta^{n-1}(\ul{\eta}) \circ \theta^{n-2}(\ul{\eta}) \circ ... \circ \theta(\ul{\eta}) \circ \ul{\eta})
$$

\vsp 
\begin{pro} 
For any two $\phi$-invariant von-Neumann measurements $\ul{\zeta}$ and $\ul{\eta}$ for $(\clm,\theta,\phi)$ we have 
\be 
H_{\phi}(\ul{\zeta}_{(n)}|\ul{\eta}_{(n)}) \le n H_{\phi}(\ul{\zeta} | \ul{\eta})
\ee
\end{pro} 

\vsp 
\begin{proof} 
We use Proposition 2.3 to verify the following inequalities:
$$H_{\phi}(\ul{\zeta}_{(n)}| \ul{\eta}_{(n)})$$
$$=H_{\phi}(\theta^{n-1}(\ul{\zeta}) ..\circ \theta(\ul{\zeta}) \circ \ul{\zeta} \circ \theta^{n-1}(\ul{\eta})\circ ... \circ \theta(\ul{\eta}) \circ \ul{\eta})$$
$$- H_{\phi}(\theta^{n-2}(\ul{\zeta}) \circ ... \circ \theta(\ul{\zeta}) \circ \ul{\zeta} 
\circ \theta^{n-1}(\ul{\eta}) \circ \theta^{n-2}(\ul{\zeta}) ... \circ \theta(\ul{\eta}) \circ \ul{\eta})$$

$$ + H_{\phi}(\theta^{n-2}(\ul{\zeta}) \circ ... \circ \theta(\ul{\zeta}) \circ \ul{\zeta} 
\circ \theta^{n-1}(\ul{\eta})\circ ..\circ \theta(\ul{\eta}) \circ \ul{\eta})
- H_{\phi}(\theta^{n-1}(\ul{\eta}) \circ \theta^{n-2}(\ul{\eta}) .. \circ \theta(\ul{\eta}) \circ \ul{\eta})$$

$$ \le H_{\phi}(\theta^{n-1}(\ul{\zeta}) \circ \theta^{n-2}(\ul{\zeta}) ... \circ \theta(\ul{\zeta}) \circ \ul{\zeta} 
| \theta^{n-1}(\ul{\eta}))$$
$$+H_{\phi}(\ul{\zeta}_{(n-1)}|\ul{\eta}_{(n-1)})$$
( by Proposition 2.3 since $\theta^{n-1}(\ul{\eta})$ commutes with $\ul{\eta}_{(n-1)}$)  

$$\le H_{\phi}(\theta^{n-1}(\ul{\zeta})|\theta^{n-1}(\ul{\eta})) + H_{\phi}(\ul{\zeta}_{(n-1)}|\ul{\eta}_{(n-1)})$$
( by Proposition 2.3, where we used the commuting property of $\theta^{n-1}(\ul{\zeta})$ with $\ul{\zeta}_{(n-1)}$
Thus we conclude (23) by induction on $n$ since $H_{\phi}(\theta^n(\ul{\zeta})|\theta^n(\ul{\eta})) = 
H_{\phi}(\ul{\zeta}|\ul{\eta})$ for all $n \ge 1$.
\end{proof} 

\vsp 
\begin{rem} 
In general inequality (23) is false even when $\clm$ is abelian. We will get back to this point with counter examples after proving our main result ( Theorem 3.11 ) of this section. 
\end{rem}

\vsp 
A partition $\ul{\eta}=(\eta_j)$ of unity in $\clm$ is called {\it projection valued von-Neumann measurement} if $\eta_j(x)=\hat{\eta}_j x \hat{\eta}_j$, where $(\hat{\eta}_j)$ is a family of orthogonal projections in $\clm$ so that $\sum_j\hat{\eta}_j=I$. In such a case, $\eta \eta_j = \eta_j$ for all $\eta_j \in \ul{\eta}$, where $\eta=\sum_i \eta_i$ and $\eta^2=\eta$ i.e. $\eta$ is a norm one projection from $\clm$ with the fixed point von-Neumann sub-algebra $\clm_{\eta}$ defined by $\clm_{\eta}=\{x \in \clm: \eta(x)= x \}$. So the fixed point sub-algebra $\clm_{\eta}$ contains the commutative von-Neumann sub-algebra $\clm_{\hat{\ul{\eta}}}$ generated by $(\hat{\eta}_j)$. 

\vsp 
Let $\ul{\eta}$ be a projection valued partition of unity in $\clm$ such that the family $\{\theta^k(\ul{\eta}): k \in \IZ \}$ of measurements are mutually commuting i.e. $\{\theta^k(\ul{\hat{\eta}})=(\theta^k(\hat{\eta}_j)): k \in \IZ,\; \eta_j \in \ul{\eta} \}$ is a commuting family of projections in $\clm$. Thus 
$$\clm_{\hat{\ul{\eta}}(\theta)} = \vee_{k \in \IZ} \clm_{\theta^k(\hat{\ul{\eta}})}$$ 
is a $\theta$ invariant commutative von-Neumann sub-algebra of $\clm$ and 
$\clm_{\hat{\ul{\eta}}(\theta)} \subseteq \{ x \in \clm: \eta(\theta^k(x))=\theta^k(x),\;\forall k \in \IZ \}.$ A projection valued von-Neumann measurement $(\ul{\eta})$ for an automorphism $\theta:\clm \raro \clm$ is called {\it maximal} if 
$$\clm_{\hat{\ul{\eta}}(\theta)} = \{ x \in \clm: \eta(\theta^k(x))=\theta^k(x),\;\forall k \in \IZ \}$$ 

\vsp 
The basic question that we need to address for a given dynamics $(\clm,\theta,\phi)$ how efficiently we can approximate quantum information gained by repeated measurements $\{\ul{\zeta}_{(k)}:k \in \IZ\}$ by projection valued von-Neumann measurements $\{\ul{\eta}_{(k)}:k \in \IZ \}$. 

\vsp 
Let $\clm$ be a von-Neumann algebra acting on a Hilbert space $\clh$ and $\phi$ be 
a faithful normal state $\phi$ of $\clm$. Let $(\clh_{\phi},\pi_{\phi},\zeta_{\phi})$ be the 
GNS space of $(\clm,\phi)$. So $\pi_{\phi}(\clm)$ is a von-Neumann algebra with a cyclic and 
separating vector $\zeta_{\phi}$ of unit length in the GNS Hilbert space $\clh_{\phi}$. We consider the algebraic tensor product $\clm \otimes \clm$ and the positive linear functional $\hat{\phi}$ 
defined by extending linearly the following map
\be 
x \otimes y \raro \hat{\phi}( x \otimes y)=\langle \zeta_{\phi}, x \clj_{\phi} y \clj_{\phi}   \rangle
\ee
The positive map $\hat{\phi}$ has a unique linear extension to the norm closure $\hat{\clm} = \overline{\clm \otimes \clm}$ of $\clm \otimes \clm$ in $\clb(\clh \otimes \clh)$ and it gives a 
state on $C^*$-algebra $\hat{\clm}$. Furthermore, for a given completely positive map $\tau:\clm \raro \clm$, with $\phi \tau \le \phi$ and $\tau(I) \le I$, there exists a unique completely positive map
$\tilde{\tau}$ on $\clm$ such that 
\be 
\langle \clj_{\phi} x \clj_{\phi} \zeta_{\phi}, \tau(y) \zeta_{\phi} \rangle  
= \langle \clj_{\phi} \tilde{\tau}(x) \clj_{\phi}, y \zeta_{\phi} \rangle
\ee
for all $x,y \in \clm$, where $\tilde{\tau}(I) \le I$ and $\phi \tilde{\tau} \le \phi$. In case $\tau$ is
a norm one projection i.e. $\tau^2 = \tau$ and $\tau(I)=I$ then $\tilde{\tau}=\tau$. For details we refer to [AC], Chapter 8 in [OP] and [Mo1].    

\vsp 
So for a $\phi$-invariant measurement $\ul{\zeta}=(\zeta_i)$, each $\zeta_i$ satisfies 
$\phi \zeta_i \le \phi \zeta = \phi$ and $\zeta_i(I) \le \zeta(I)=I$. Thus $\tilde{\zeta}=(\tilde{\zeta}_i)$ is also a $\phi$-invariant measurement in $\clm$.  
Let $\ul{\eta}=(\ul{\eta}_j)$ be a projection valued measurement in $\clm$ such that 
$\phi \eta = \phi$. Then $\eta$ is a $\phi$-invariant norm one projection onto the fixed 
point algebra $\clm_{\eta}$ and $\eta$ commutes with the modular group $\sigma^{\phi}=(\sigma^{\phi}_t)$. Thus the family $\hat{\eta}_j$ of projections in $\clm$ are also $(\sigma^{\phi}_t)$-invariant. 
In particular, we have $\tilde{\eta}_j = \eta_j$ for all $j$. 

\vsp 
Furthermore, $\IE_{\hat{\ul{\eta}}}:\clm \raro \clm_{\hat{\ul{\eta}}}$ is a $\phi$ invariant norm one projection defined by 
\be 
\IE_{\hat{\ul{\eta}}}(x)= \sum_j { \phi(\hat{\eta}_j x \hat{\eta}_j) \over \phi(\hat{\eta}_j) } \hat{\eta}_j 
\ee

\vsp 
\begin{pro} 
Let $\ul{\eta}=(\eta_j)$ be a $\phi$-invariant projection valued measurements in $\clm$ then 
for any $\phi$-invariant measurement $\ul{\zeta}=(\zeta_i)$, we have  
\be 
H_{\phi}(\ul{\zeta} \circ \ul{\eta}) - H_{\phi}(\ul{\eta}) = - \sum_i \phi(\IE_{\hat{\ul{\eta}}}(\tilde{\zeta}_i(1)) ln \IE_{\hat{\ul{\eta}}}(\tilde{\zeta}_i(I)))  +  \sum_i S( \phi \zeta_i | \clm_{\eta},\phi |\clm_{\eta})
\ee 
\end{pro} 

\vsp 
\begin{proof} 
\vsp 
We compute as in Proposition 2.2 that
$$H^q_{\phi}(\ul{\zeta} \circ \ul{\eta}) - H^q_{\phi}(\ul{\eta}) $$
$$ = \sum_{i,j} S(\phi \zeta_i \eta_j, \phi \eta_j)$$
\be 
= \sum_{i,j} S(\phi \zeta_i \eta \eta_j, \phi \eta \eta_j)
\ee
(since $\eta_j \eta= \eta_j$)
\be 
\le \sum_i S(\phi \zeta_i \eta, \phi \eta)
\ee
by monotonicity property of relative entropy. On the other hand, join convexity and scaling property says for any finite family of normal states $\omega_i$ and $\phi_i$, we have 
\be 
S(\sum_i \omega_i, \sum_i \phi_i) \le \sum_i S(\omega_i,\phi_i)
\ee
and thus 
$$H^q_{\phi}(\ul{\zeta} \circ \ul{\eta}) - H^q_{\phi}(\ul{\eta}) $$
$$\ge \sum_i S(\phi \zeta_i \eta, \phi \eta)$$
$$=\sum_i S(\phi \zeta_i | \clm_{\eta}, \phi |\clm_{\eta})$$ 
Thus we have the following equality:
\be 
H^q(\ul{\zeta} \circ \ul{\eta})-H^q(\ul{\eta}) = \sum_i S( \phi \zeta_i | \clm_{\eta},\phi |\clm_{\eta})
\ee

\vsp 
We also compute the following identities:
$$H^c_{\phi}(\ul{\eta} \circ \ul{\zeta}) -H^c_{\phi}(\ul{\eta})$$
$$= \sum_{i,j} -\phi(\eta_j \zeta_i(I))ln \phi(\eta_j \zeta_i(I)) + \phi(\eta_j(I)ln(\phi(\eta_j(I)$$
$$= \sum_{i,j} - \phi( \hat{\eta}_j \eta \circ \zeta_i(I) \hat{\eta}_j) ln \phi (\hat{\eta}_j \eta \circ \zeta_i(I) \hat{\eta}_j) +  \phi(\hat{\eta}_j) ln(\phi(\hat{\eta}_j)$$
$$ = \sum_{i,j} - \phi( \hat{\eta}_j \zeta^{\eta}_i(I) \hat{\eta}_j) ln \phi(\hat{\eta}_j \zeta^{\eta}_i(I) \hat{\eta}_j) + 
\phi(\hat{\eta}_j) ln(\phi(\hat{\eta}_j)$$
$$= \sum_{i,j} \phi({\eta}_j) {\phi ( \hat{\eta}_j \zeta^{\eta}_i(I) \hat{\eta}_j )  \over \phi(\hat{\eta}_j) } ln \phi({ \hat{\eta}_j \zeta^{\eta}_i(I) \hat{\eta}_j \over \phi(\tilde{\eta}_j)})$$
$$ = \sum_i - \phi \IE_{\hat{\ul{\eta}}}[\zeta_i^{\eta}(I)] ln (\IE_{\hat{\ul{\eta}}}[\zeta_i^{\eta}(I)])$$

\vsp 
$$H^c_{\phi}(\ul{\zeta} \circ \ul{\eta})-H^c_{\phi}(\ul{\eta})$$
$$=H^c_{\phi}(\ul{\tilde{\eta}} \circ \ul{\tilde{\zeta}}) - H^c_{\phi}(\ul{\tilde{\eta}})$$
\be 
= \sum_i - \phi \IE_{\hat{\ul{\eta}}}[\tilde{\zeta}_i^{\eta}(I)] ln (\IE_{\hat{\ul{\eta}}}[\tilde{\zeta}_i^{\eta}(I)])
\ee
where we have used $\eta_j=\tilde{\eta}_j$. 

\vsp 
We combine equalities given in (31) and (32) to complete the proof for equality given in (27). 
\end{proof} 

\vsp 
Let the triplet $(\clm,\theta,\phi)$ be as in Proposition 3.2 and $\ul{\eta}$ be a $\phi$-invariant projection valued partition of unity in $\clm$ as in Proposition 3.6. Furthermore, let $\phi_0$ be possibly an another faithful normal state on $\clm$ such that $\phi_0 \theta = \phi_0$ and $\phi_0 \eta=\phi_0$. Let $\cln_{\hat{\ul{\eta}}}$ be the maximal sub-algebra of $\clm$ for which  
\be 
\clm_{\hat{\ul{\eta}}}' \bigcap \cln_{\hat{\ul{\eta}}} \subseteq \clm_{\hat{\ul{\eta}}},
\ee
and $\cln_{\hat{\ul{\eta}}}$ is preserved by the modular automorphism $(\sigma^{\phi_0}_t)$ associated with $\phi_0$ on $\clm$. 
That such a maximal von-Neumann algebra exists follows trivially once we verify that $\sigma_t^{\phi_0}(\clm_{\hat{\ul{\eta}}})=\clm_{\hat{\ul{\eta}}}$ for all $t \in \IR$. Since $\phi_0 \eta=\phi_0$ and $\eta$ is a norm one projection, $\eta$ commutes with the modular group $(\sigma_t^{\phi_0})$ by a Theorem of M. Takasaki [Ta1], we verify that $\{\sigma_t^{\phi_0}(\hat{\eta}_j) \}'= \{ x \in \clm : \eta(x)= x \}$ for each $t \in \IR$. In fact we have $\sigma_t^{\phi_0}(\hat{\eta}_j)=\hat{\eta}_j$ for all $t \in \IR$ and $\eta_j$. So by our construction, we have $\clm_{\hat{\ul{\eta}}} \subseteq \cln_{\hat{\ul{\eta}}}$. By the same theorem of M. Takasaki [Ta1], there exists a $\phi_0$ invariant norm one projection $\IE_{\cln_{\ul{\hat{\eta}}}}$ onto $\cln_{\ul{\hat{\eta}}}$ such that $\phi_0 \IE_{\cln_{\ul{\hat{\eta}}}}=\phi_0$ on $\clm$. 

\vsp 
Thus we can choose recursively maximal sub-algebras $(\cln_{\hat{\ul{\eta}}_{(n)} }:n \ge 1)$ of $\clm$
satisfying for each $n \ge 1:$
$$\clm'_{\hat{\ul{\eta}}_{(n+1)}} \bigcap \cln_{\hat{\ul{\eta}}_{(n)}} 
\subseteq \clm'_{\hat{\ul{\eta}}_{(n)}} \bigcap \cln_{\hat{\ul{\eta}}_{(n)}} \subseteq \clm_{\hat{\ul{\eta}}_{(n)}} \subseteq \clm_{\hat{\ul{\eta}}_{(n+1)}}$$
\be 
\clm_{\hat{\ul{\eta}}_{(n)}}' \bigcap \cln_{\hat{\ul{\eta}}_{(n)}} \subseteq \clm_{\hat{\ul{\eta}}_{(n)}}
\ee
and 
$$\sigma^{\phi_0}_t(\cln_{\hat{\ul{\eta}}_{(n)}})=\cln_{\hat{\ul{\eta}}_{(n)}},\; t \in \IR$$
such that 
$$\cln_{\hat{\ul{\eta}}_{(n)}} \subseteq \cln_{\hat{\ul{\eta}}_{(n+1)}}$$ 
We have 
$$\cln_{\hat{\ul{\eta}}(\theta)^+} = \vee_{0 \le n < \infty} \cln_{\hat{\ul{\eta}}_{(n)}} \subseteq \clm$$ 
and $\sigma_t^{\phi_0}(\cln_{\ul{\hat{\eta}}(\theta)^+})=\cln_{\ul{\hat{\eta}}(\theta)^+}$. 
                                                                                                                                                                                                                                                                                                                                                                                                                                                                                                                                                                                                                                                                                                                                                                                                                                                                                                                                                                                                                                                                                                                                                                                                                                                                                                                                                                                                                                                                                                                                                                                                                                                                                                                                                                                                                                                                                                                                                                                                                                                                                                                                                                                                                                                                                                                                                                                                                                                                                                                                                                                                                                                                                                                                                                                                                                                                                                                                                                                                                                                                                                                                                                                                                                                                                                                                                                                                                                                                                                                                                                                                    
\vsp 
More generally, we can choose recursively $(\sigma_t^{\phi_0})$-invariant maximal sub-algebras $\cln_{\hat{\ul{\eta}}^{(m)}_{(n)}}:-\infty < m \le n < \infty$ of $\clm$ such that 
$$\clm_{\hat{\ul{\eta}}^{(m)}_{(n)}}' \bigcap \cln_{\hat{\ul{\eta}}^{(m)}_{(n)}} \subseteq \clm_{\hat{\ul{\eta}}^{(m)}_{(n)}}$$
and 
$$\cln_{\hat{\ul{\eta}}^{(m)}_{(n)}} \subseteq \cln_{\hat{\ul{\eta}}^{(m')}_{(n')}}$$ 
if $m' \le m \le n \le n'$, where $\ul{\eta}^{(m)}_{(n)} = \theta^m(\ul{\eta})\circ ...\circ \theta^n(\ul{\eta})$ for $-\infty < m \le n < \infty$. It is evident that 
$$\cln_{\hat{\ul{\eta}}(\theta)} = \vee_{- \infty < m \le n < \infty} \cln_{\hat{\ul{\eta}}^{(m)}_{(n)}} \subseteq \clm$$ 
and $\sigma_t^{\phi_0}(\cln_{\hat{\ul{\eta}}(\theta)})=\cln_{\hat{\ul{\eta}}(\theta)}$. 
                                                                                                                                                                                                                                                                                                                                                                                                                                                                                                                                                                                                                                                                                                                                                                                                                                                                                                                                                                                                                                                                                                                                                                                                                                                                                                                                                                                                                                                                                                                                                                                                                                                                                                                                                                                                                                                                                                                                                                                                                                                                                                                                                                                                                                                                                                                                                                                                                                                                                                                                                                                                                                                                                                                                                                                                                                                                                                                                                                                                                                                                                                                                                                                                                                                                                                                                                                                                                                                                                                                                                                                                    
                                                                                                                                                                                                                                                                                                                                                                                                                                                                                                                                                                                                                                                                                                                                                                                                                                                                                                                                                                                                                                                                                                                                                                                                                                                                                                                                                                                                                                                                                                                                                                                                                                                                                                                                                                                                                                                                                                                                                                                                                                                                                                                                                                                                                                                                                                                                                                                                                                                                                                                                                                                                                                                                                                                                                                                                                                                                                                                                                                                                                                                                                                                                                                                                                                                                                                                                                                                                                                                                                                                                                                                                                                                                                                                                                                                                                                                                                                                                                                                                                                                                                                                                                                                                                                                                                                                                                                                                                                                                                                                                                                                                                                                                                                                                                                                                                                                                                                                                                                                                                                                                                                                                                                                                                                                                                                                                                                                                                                                                                                                                                                                                                                                                                                                                                                                                                                                                                                              \vsp 
Since $\IE_{\cln_{\hat{\ul{\eta}}_{(n+m)}}} \IE_{\cln_{\hat{\ul{\eta}}_{(m)}}} = \IE_{\cln_{\hat{\ul{\eta}}_{(n)}}}$, by Proposition 2.1 (h) we have the following sub-additive property: 
$$S_{n+m} \le S_n + S_m$$ 
for 
$$S_n= S(\hat{\phi}|\cln_{\hat{\ul{\eta}}_{(n)}} \otimes \cln_{\hat{\ul{\eta}}_{(n)}}, \hat{\phi} \eta_{(n)} \otimes \eta_{(n)}|\cln_{\hat{\eta}_{(n)}} \otimes \cln_{\hat{\eta}_{(n)}})$$
We define mean relative entropy $s(\hat{\phi},\hat{\phi} \circ \eta \otimes \eta)$ by 
$$s(\hat{\phi},\hat{\phi} \circ \eta \otimes \eta) = \mbox{lim}_{n \raro \infty}{1 \over n} S_n$$

\vsp 
By monotonicity of relative entropy, we also have 
$$\sum_i S(\phi \zeta_i|\cln_{\hat{\ul{\eta}}}, \phi \zeta_i \eta|\cln_{\hat{\ul{\eta}}})$$
$$\le \sum_i S(\hat{\phi} \eta \zeta_i \otimes I| \cln_{\hat{\ul{\eta}}} \otimes \cln_{\hat{\ul{\eta}}}, \hat{\phi} \eta \zeta_i \eta  \otimes I)|\cln_{\hat{\ul{\eta}}} \otimes \cln_{\hat{\ul{\eta}}})$$
$$ = \sum_i S(\hat{\phi} I \otimes \tilde{\zeta}_i \eta |\cln_{\hat{\ul{\eta}}} \otimes \cln_{\hat{\ul{\eta}}}, \hat{\phi} \eta \otimes \tilde{\zeta}_i \eta|\cln_{\hat{\ul{\eta}}} \otimes \cln_{\hat{\ul{\eta}}})$$
$$ \le S(\hat{\phi}|\cln_{\hat{\ul{\eta}}} \otimes \cln_{\hat{\ul{\eta}}}, \hat{\phi} \eta \otimes I |\cln_{\hat{\ul{\eta}}} \otimes \cln_{\hat{\ul{\eta}}})$$
\be 
\le S(\hat{\phi}|\cln_{\hat{\ul{\eta}}} \otimes \cln_{\hat{\ul{\eta}}}, \hat{\phi} \eta \otimes \eta |\cln_{\hat{\ul{\eta}}} \otimes \cln_{\hat{\ul{\eta}}})
\ee
In the last line we have used 
                                                                                                                                                                                                                                                                                                                                                                                                                                                                                                                                                                                                                                                                                                                                                                                                                                                                                                                                                                                                                                                                                                                                                                                                                                                                                                                                                                                                                                                                                                                                                                                                                                                                                                                                                                                                                                                                                                                                                                                                                                                                                                                                                                                                                                                                                                                                                                                                                                                                                                                                                                                                                                                                                                                                                                                                                                                                                                                                                                                                                                                                                                                                                                                                                                                                                                                                                                                                                                                                                                                                                                                                                                                                                                                                                                                                                                                                                                                                                                                                                                                      $\eta^2=\eta$ and $\hat{\phi} \eta \otimes I = \hat{\phi} I \otimes \eta$.

\vsp 
\begin{pro} 
Let $\ul{\zeta}$ and $\ul{\eta}$ be measurements in $\clm$ as in Proposition 3.6. Then  
\be 
0 \le H_{\phi}(\ul{\zeta})-H_{\phi}(\ul{\zeta} \circ \eta) = \sum_i S(\phi \zeta_i, \phi \zeta_i \eta ) 
\ee 
and 
\be 
0 \le h_{\phi}(\theta,\ul{\zeta}) - \mbox{lim}_{n \raro \infty} {1 \over n} H_{\phi}(\ul{\zeta}_{(n)} \circ \eta_{(n)}) \le s(\hat{\phi}, \hat{\phi} \circ \eta \otimes \eta) + S(\hat{\phi},\hat{\phi} \IE_{\cln_{\hat{\ul{\eta}}(\theta)^+}} \otimes \IE_{\hat{\ul{\eta}}(\theta)^+})
\ee
\end{pro} 

\vsp 
\begin{proof} 
By Proposition 2.1 (g), we have 
$$H_{\phi}(\ul{\zeta})= \sum_i \phi(\zeta_i(I))S( {1 \over \phi \zeta_i(I) } \phi \zeta_i, \phi)$$
$$= \sum_i \phi(\zeta_i(I)) \{ S( {1 \over \phi \zeta_i(I) } \phi \zeta_i|\clm_{\eta}, \phi|\clm_{\eta}) + 
S({ 1 \over \phi \zeta_i(I)} \phi \zeta_i, { 1 \over \phi \zeta_i(I)} \phi \zeta_i \eta) \}$$
So by scaling property given in Proposition 2.1 (a), we have 
$$H_{\phi}(\ul{\zeta})= -\sum_i \phi(\zeta_i(I))ln \phi(\zeta_i(I)) + \sum_i S(\phi \zeta_i|\clm_{\eta}, \phi|\clm_{\eta}) + S(\phi \zeta_i, \phi \zeta_i \eta)$$
Thus 
$$H^q_{\phi}(\ul{\zeta}) = \sum_i S(\phi \zeta_i , \phi) $$
$$= \sum_i S(\phi \zeta_i|\clm_{\eta}, \phi|\clm_{\eta}) + S(\phi \zeta_i, \phi \zeta_i \eta)$$ 
and so 
$$H_{\phi}(\ul{\zeta})-H_{\phi}(\ul{\zeta} \circ \eta)$$
$$=H^q_{\phi}(\ul{\zeta})-H^q_{\phi}(\ul{\zeta} \circ \eta)$$
$$=\sum_i S(\phi \zeta_i, \phi)-S(\phi \zeta_i \eta, \phi)$$
$$=\sum_i S(\phi \zeta_i, \phi \zeta_i \eta)$$
$$=\sum_i S(\phi \zeta_i, \phi_0)-S(\phi \zeta_i \eta, \phi_0)$$
(by Proposition 2.1 (h) since $\phi_0 \eta = \phi_0$ as well)

$$=\sum_i S(\phi \zeta_i|\cln_{\hat{\eta}}, \phi_0|\cln_{\hat{\eta}}) + S(\phi \zeta_i, \phi \zeta_i \IE_{\cln_{\hat{\ul{\eta}}}})$$
$$- \sum_i S(\phi \zeta_i \eta |\cln_{\hat{\ul{\eta}}}, \phi_0|\cln_{\hat{\ul{\eta}}}) + S(\phi \zeta_i \eta, \phi \zeta_i \eta \IE_{\cln_{\hat{\ul{\eta}}}})$$
(by Proposition 2.1 (h) since $\phi_0 \IE_{\cln_{\hat{\ul{\eta}}}}=\phi_0$)
$$=\sum_i S(\phi \zeta_i|\cln_{\hat{\eta}}, \phi_0|\cln_{\hat{\eta}}) - S(\phi \zeta_i \eta |\cln_{\hat{\ul{\eta}}}, \phi_0|\cln_{\hat{\ul{\eta}}})$$
$$+ \sum_i S(\phi \zeta_i, \phi \zeta_i \IE_{\cln_{\hat{\ul{\eta}}}}) - S(\phi \zeta_i \eta, \phi \zeta_i \eta \IE_{\cln_{\hat{\ul{\eta}}}})$$
$$\le \sum_i S(\phi \zeta_i |\cln_{\hat{\ul{\eta}}},\phi \zeta_i \eta |\cln_{\hat{\ul{\eta}}}) + S(\phi \zeta_i, \phi \zeta_i \IE_{\cln_{\hat{\ul{\eta}}}})$$
since $$\sum_i S(\phi \zeta_i \eta, \phi \zeta_i \eta \IE_{\cln_{\hat{\ul{\eta}}}}) \ge S( \sum_i \phi \zeta_i \eta, \sum_i \phi \zeta_i \eta \IE_{\cln_{\hat{\ul{\eta}}}}) = 0 $$

\vsp 
So we have  
$$H_{\phi}(\ul{\zeta}_{(n)}) - H_{\phi}(\ul{\zeta}_{(n)} \circ \eta_{(n)})$$
$$ \le S_n + S(\hat{\phi}, \hat{\phi} \IE_{\cln_{\hat{\eta}_{(n)}}} \otimes \IE_{\cln_{\hat{\eta}_{(n)}}} ) $$
Since  $S(\hat{\phi}, \hat{\phi} \IE_{\cln_{\hat{\eta}_{(n)}}} \otimes \IE_{\cln_{\hat{\eta}_{(n)}}} ) \downarrow 
S(\hat{\phi},\hat{\phi} \IE_{\cln_{\ul{\hat{\eta}}(\theta)^+}} \otimes \IE_{\cln_{\hat{\ul{\eta}}(\theta)^+}})$ as $n \raro \infty$. 
\end{proof} 

\vsp 
\begin{pro} 
Let $\ul{\zeta}$ and $\ul{\eta}$ be two $\phi$-invariant measurements for $(\clm,\theta,\phi)$. If $\ul{\eta}$ is a von-Neumann measurement for $(\clm,\theta,\phi)$ and  
\be 
s(\hat{\phi},\hat{\phi} \eta \otimes \eta)  < \infty 
\ee
then we have
\be 
h_{\phi}(\theta,\ul{\zeta}) - h_{\phi}(\theta,\ul{\eta}) \le H_{\phi}(\ul{\zeta}|
\ul{\eta}(\theta)) + s(\hat{\phi}, \hat{\phi} \eta \otimes \eta) + S(\hat{\phi}.\hat{\phi} \IE_{\cln_{\hat{\ul{\eta}}(\theta)}} \otimes \IE_{\cln_{\hat{\ul{\eta}}(\theta)}})
\ee 
where 
$$
H_{\phi}(\ul{\zeta}|\ul{\eta}(\theta)) = - \sum_i \phi(\IE_{\hat{\eta}(\theta)}(\tilde{\zeta}_i(1)) ln \IE_{\hat{\eta}(\theta)}(\tilde{\zeta}_i(I)))  +  \sum_i S( \phi \zeta_i | \clm_{\eta(\theta)},\phi |\clm_{\eta(\theta)})
$$
\end{pro} 

\vsp 
\begin{proof} 
\vsp 
By Proposition 3.6 we have 
$$H_{\phi}(\ul{\zeta}_{(n)} \circ \eta_n) - H_{\phi}(\ul{\eta}_n)$$
$$\le H_{\phi}(\ul{\zeta}_{(n)}|\ul{\eta}_{(n)})$$
$$\le n H_{\phi}(\ul{\zeta}|\ul{\eta})$$  

\vsp 
However by Proposition 3.7, we have 
$$H_{\phi}(\ul{\zeta}_{(n)})-H_{\phi}(\ul{\zeta}_{(n)} \circ \eta_{(n)})$$
$$ \le S_n +  S(\hat{\phi}, \hat{\phi} \IE_{\cln_{\hat{\eta}_{(n)}}} \otimes \IE_{\cln_{\hat{\eta}_{(n)}}})$$
Thus we arrive at 
$$h_{\phi}(\theta,\ul{\zeta})-h_{\phi}(\theta,\ul{\eta}) \le H_{\phi}(\ul{\zeta}|\ul{\eta}) + s(\hat{\phi}, \hat{\phi} \eta \otimes \eta) + S(\hat{\phi},\hat{\phi} \IE_{\cln_{\hat{\ul{\eta}}(\theta)^+}} \otimes \IE_{\cln_{\hat{\ul{\eta}}(\theta)^+}})$$

\vsp 
Since $h_{\phi}(\theta,\ul{\eta})=h_{\phi}(\theta, \ul{\eta}^{(m)}_{(n)})$ for $\ul{\eta}^{(m)}_{(n)}=\theta^{n-1}(\ul{\eta}) \circ ...\ul{\eta} \circ \theta^{-1}(\ul{\eta}) ..\circ \theta^m(\ul{\eta})$, by the same line of argument and $s(\hat{\phi},\hat{\phi} \eta \otimes \eta)=s(\hat{\phi},\hat{\phi} \eta^{(m)}_{(n)} \otimes \eta^{(m)}_{(n)})$ for all $-\infty < m < n < \infty$, we get 
$$h_{\phi}(\theta,\ul{\zeta})-h_{\phi}(\theta,\ul{\eta}) \le H_{\phi}(\ul{\zeta}|\ul{\eta}^{(m)}_{(n)}) + s(\hat{\phi},\hat{\phi} \eta \otimes \eta)+ S(\hat{\phi},\hat{\phi} \IE_{\cln_{\theta^m(\eta(\theta)^+)}} \otimes \IE_{\cln_{\theta^m(\eta(\theta)^+)}})$$
Taking limit $n \raro \infty$ and $m \raro -\infty$, we conclude a proof for (39).  
\end{proof} 

\vsp 
\begin{pro} 
For any $\phi$-invariant von-Neumann measurement $\ul{\zeta}$ and maximal von-Neumann measurement $\ul{\eta}$ for dynamics $(\clm,\theta,\phi)$ we have 
\be 
H_{\phi}(\ul{\zeta}|\ul{\eta}(\theta))=0
\ee 
\end{pro} 

\vsp 
\begin{proof} 
We will use Proposition 3.6. The commutative von-Neumann sub-algebra $\clm_{\hat{\ul{\eta}}(\theta)}$ being maximal in $\clm$ i.e. $\clm_{\eta(\theta)}=\clm_{\hat{\ul{\eta}}(\theta)}$, we have $\IE_{\hat{\eta}(\theta)}=\eta(\theta)$ and so 
$$S(\phi \zeta_i |\clm_{\hat{\ul{\eta}}(\theta)},\phi|\clm_{\hat{\ul{\eta}}(\theta)})$$
$$=\phi(\IE_{\hat{\eta}(\theta)}(\tilde{\zeta}_i(I))ln \IE_{\hat{\eta}(\theta)}(\tilde{\zeta}_i(I)))$$ 
for each $\zeta_i$. Thus (27) gives the required result. 
\end{proof} 

\vsp 
\begin{pro} 
Let $(\clm,\theta,\phi)$ be a $W^*$-dynamical system with a faithful normal invariant state for the automorphism $\theta$. If there exist a maximal von-Neumann partition of unity in $\clm$ for
the dynamics $(\clm,\theta,\phi)$ with $s(\hat{\phi}, \hat{\phi} \eta \otimes \eta) < \infty$ and $\cln_{\hat{\ul{\eta}}(\theta)}=\clm$ then 
$$h_{\phi}(\theta,\ul{\zeta}) \le h_{\phi}(\theta,\ul{\eta})+  s(\hat{\phi}, \hat{\phi} \eta \otimes \eta)$$
for any $\phi$-invariant von-Neumann measurement $\ul{\zeta}$ in $\clm$ 
with $H_{\phi}(\ul{\zeta}) < \infty$.  
\end{pro} 

\vsp 
A $\phi$-invariant measurement $\ul{\zeta}$ is called {\it admissible} for the dynamics 
$(\clm,\theta,\phi)$ if there exists an integer $k \ge 1$ such $\ul{\zeta}$ is a von-Neumann 
measurement for $(\clm,\theta^k,\phi)$. We note that $\ul{\zeta} \circ \ul{\eta}$ is also an admissible measurement for $(\clm,\theta,\phi)$ if $\ul{\zeta}$ and $\ul{\eta}$ are so i.e. they form an algebra 
under composition $\circ$. For a given $\phi$-invariant admissible measurement $\ul{\zeta}$ for $(\clm,\theta)$, 
we set 
\be 
h_{\phi}(\theta,\ul{\zeta})= \mbox{sup} {1 \over k } h_{\phi}(\theta^k,\ul{\zeta})
\ee
where sup is taken over all possible values for $k$ for which $\ul{\zeta}$ are von Neumann measurements for $(\clm,\theta^k,\phi)$. 

\vsp 
\begin{thm} 
Let $(\clm,\theta,\phi)$ be a $W^*$-dynamical system with a faithful normal invariant state for the automorphism $\theta$. If there exist a maximal von-Neumann partition of unity in $\clm$ for
the dynamics $(\clm,\theta,\phi)$ such that $\cln_{\hat{\ul{\eta}}(\theta)}=\clm$ for some faithful normal state $\omega_0$ then 
\be 
h_{\phi}(\theta,\ul{\zeta}) \le h_{\phi}(\theta,\ul{\eta}) +s(\hat{\phi}, \hat{\phi} \eta \otimes \eta)
\ee
for any $\phi$-invariant admissible measurement $\ul{\zeta}$ for $(\clm,\theta,\phi)$  
with $H_{\phi}(\ul{\zeta}) < \infty$.  
\end{thm} 

\vsp 
If we intend to obtain information of a dynamical system $(\clm,\theta,\phi)$ by repeated measurement, then the set of $\phi$-invariant admissible measurements $\clp_{\phi,a}(\theta)$ for $(\clm,\theta,\phi)$ 
is a natural class. We set 
\be 
h_{\phi}(\theta) = \mbox{sup}_{\ul{\zeta} \in \clp_{\phi,a}(\theta) }  h_{\phi}(\theta, \ul{\zeta})
\ee 
The value $h_{\phi}(\theta)$ could be infinite. As in the classical theory, we expect this notion to play an important role when $0 < h_{\phi}(\theta) < \infty$. In particular we have the following:

\vsp 
\begin{thm} 
The class of admissible partitions for $(\clm,\theta,\phi)$ is an invariance under isomorphism and $h_{\phi}(\theta)$ is a non-negative invariance for $(\clm,\theta,\phi)$. Furthermore we have 

\vsp 
\NI (a) $h_{\phi}(\theta^m) = m h_{\phi}(\theta);$

\NI (b) If $\theta^m=I$ for some $m \in \IZ$, then $h_{\phi}(\theta)=0$ or $\infty$. 
\end{thm} 

\vsp 
\begin{proof} 
Let $\ul{\zeta}$ be an admissible partition for $(\clm,\theta^m,\omega)$. So there exists a $n \ge 1$ 
so that $\ul{\zeta}$ is von-Neumann partition for $(\clm,\theta^{mn},\omega)$ i.e. $\{ \theta^{nmk}(\ul{\zeta}): k \in \IZ \}$ are mutually commuting. Thus $\theta^n(\ul{\zeta})$ is a von-Neumann partition for $(\clm,\theta^m,\phi)$ and so $(\theta^n(\ul{\zeta}))$ is an admissible partition for $(\clm,\theta,\phi)$. By (41), we have  
$$h_{\phi}(\theta,\theta^n(\ul{\zeta})) \ge {1 \over m} h_{\phi}(\theta^m,\theta^n(\ul{\zeta}))$$
Since $h_{\phi}(\theta^m,\theta^n(\ul{\zeta})=h_{\phi}(\theta^m,\ul{\zeta})$ for any $n$, we have 
$$h_{\phi}(\theta^m, \ul{\zeta}) \le m h_{\phi}(\theta,\theta^n(\ul{\zeta})) \le m h_{\phi}(\theta)$$ 
for all admissible partition $\ul{\zeta}$ for $(\clm,\theta^m,\phi)$. This shows $h_{\phi}(\theta^m) 
\le m h_{\phi}(\theta)$. 
 
\vsp 
Conversely, let $\ul{\zeta}$ be an admissible partition for $(\clm,\theta,\phi)$. So there exists a 
$n \ge 1$ so that $\ul{\zeta}$ is a von-Neumann partition for $(\clm,\theta^n,\phi)$ for some $n \ge 1$ i.e. 
$\{ \theta^{kn}(\ul{\zeta}): k \in \IZ \}$ are mutually commuting. Now the partition 
$\theta^{(m-1)n}(\ul{\zeta}) \circ...\theta^n(\ul{\zeta}) \circ \ul{\zeta}$ is a von-Neumann partition for $(\clm,\theta^{mn},\phi)$ and thus $\theta^{(m-1)n}(\ul{\zeta}) \circ...\theta^n(\ul{\zeta}) \circ \ul{\zeta}$ is an admissible partition for $(\clm,\theta^m,\phi)$ and 
$$h_{\phi}(\theta^m)$$
$$\ge h_{\phi}(\theta^m, \theta^{(m-1)n}(\ul{\zeta}) \circ...\theta^n(\ul{\zeta}) \circ \ul{\zeta})$$  
$$ \ge {1 \over n } h_{\phi}(\theta^{mn}, \theta^{(m-1)n}(\ul{\zeta}) \circ...\theta^n(\ul{\zeta}) \circ \ul{\zeta})$$  
$$\ge m h_{\phi}(\theta,\ul{\zeta}) $$  
Since the inequality holds for any admissible partition for $(\clm,\theta,\phi)$, we get $h_{\phi}(\theta^m) \ge m h_{\phi}(\theta)$. This shows equality in (a)

\vsp 
The statement (b) follows trivially once we show $h_{\phi}(I)=0$ or $\infty$. Suppose $h_{\phi}(I) < \infty$. Then by (a), we have 
$n h_{\phi}(I)=h_{\phi}(I)$ for all $n \ge 1$ and so $h_{\phi}(I)=0$.
\end{proof}

\vsp 
\begin{rem}
Theorem 3.11 says that dynamical entropy $h_{\phi}(\theta)$ of $(\clm,\theta,\phi)$ satisfies 
\be 
h_{\phi}(\theta,\ul{\eta}) \le h_{\phi}(\theta) \le h_{\phi}(\theta,\ul{\eta}) + s(\hat{\phi},\hat{\phi} \eta \otimes \eta) 
\ee
provided there exists a maximal projection valued von-Neumann partition $\ul{\eta}$ of unity for $(\clm,\theta,\phi)$ with 
$\cln_{\hat{\ul{\eta}}(\theta)}=\clm$. If the maximal von-Neumann partition $\ul{\eta}$ for $(\clm,\theta,\phi)$ is also a Kolmogorov partition 
for the maximal abelian von-Neumann algebra $\clm_{\hat{\ul{\eta}}(\theta)}$ in the classical sense, then $h_{\phi}(\theta,\ul{\eta}) > 0$. 
Furthermore, if $\clm$ is a commutative von-Neumann algebra, then we have $\eta=I$ and so $s(\hat{\phi}, \hat{\phi} \eta \otimes \eta)=0$
Thus Theorem 3.11 clearly shows that this notion is indeed a generalization of Kolmogorov-Sinai theorem [Pa]. 
\end{rem}

\vsp 
Given a $C^*$-dynamical system $(\cla,\theta,\omega)$, we consider the $W^*$ dynamical system $(\bar{\cla},\bar{\theta},\bar{\omega})$, where $\bar{\cla}=\cla^{**}$ is double dual of $\cla$ i.e. universal von-Neumann algebra over $\cla$ and $\bar{\omega}$ be the unique normal extension of $\cla$ to $\bar{\cla}$ invariant for the induced auto-morphism $\bar{\theta}$
on $\bar{\cla}$ given by $\bar{\theta}(\bar{\pi}(x))=\bar{\pi}(\theta(x))$ for all $x \in \cla$, where $\bar{\pi}:\cla \raro \bar{\cla}$ is the universal representation. Two $C^*$-dynamical systems $(\cla_1,\theta_1,\omega_1)$ and $(\cla_2,\theta_2,\omega_2)$ are called {\it isomorphic} if there exists an automorphism $\alpha: \bar{\cla}_1 \raro \bar{\cla}_2$ such that $\alpha \bar{\theta}_1 = \bar{\theta}_2 \alpha$ and $\bar{\omega}_2 \alpha=\omega_1$. We define dynamical entropy $h_{\omega}(\theta)$ of
a $C^*$-dynamics $(\cla,\theta,\omega)$ by equating to the dynamical entropy $h_{\bar{\omega}}(\bar{\theta})$ of the $W^*$ dynamical system $(\bar{\cla},\bar{\theta},\bar{\omega})$. It is obvious that $h_{\omega}(\theta)$ is an invariance for the dynamics $(\cla,\theta,\omega)$. It is also evident by our definition $h_{\omega}(\theta)=h_{\phi_{\omega}}(\theta_{\omega})$, where $\theta_{\omega}$ is the automorphism on $\pi_{\omega}(\cla)''$ defined by extending the map $\pi_{\omega}(x) \raro \pi_{\omega}(\theta(x))$ for all $x \in \cla$ with invariant normal state $\phi_{\omega}(X)=\langle \zeta_{\omega}, X \zeta_{\omega} \rangle$
for all $X \in \pi_{\omega}(\cla)''$. 

\vsp
\section{ Abelian partition and Kolmogorov Sinai theorem revisited:}

\vsp
Let $(\Omega,\clf,\mu)$ be a probability space and $\theta$ be a one-to-one and onto measurable map on 
$\Omega$ so that $\mu \circ \theta =\mu$. In this section we will characterize the maximal class of measurements 
that commute with the class of measurements associated with measurable partitions. We will show although 
the class is larger then that studied by Kolmogorov-Sinai, the dynamical entropy is same.

\vsp
\begin{pro} 
Let $\ul{\zeta}=(\zeta_i)$ be a family of positive maps so that $\sum_i\zeta_i(1)=1$ and 
$\zeta_i\eta_j=\eta_j\zeta_i$ for any $\eta_j(\psi)=\hat{\eta}_j\psi$, where $(\hat{\eta}_j)$ is a partition of the probability 
space $\Omega$ into measurable sets. Then there exists a family of bounded measurable functions $\hat{\zeta}_i$ so that 
$\zeta_i(\psi)=\hat{\zeta_i}^*\psi\hat{\zeta}_i$ and $\sum_i\hat{\zeta}_i^*\hat{\zeta}_i=1$.
\end{pro} 

\vsp
\begin{proof} 
Since $\mu(\zeta_i(\psi)) \le \mu(\psi)$, there exists a family of bounded measurable functions
$(\hat{\zeta}_i)$ so that $\mu(\zeta_i(\psi))=\mu(\hat{\zeta}^*_i\psi\hat{\zeta}_i)$ for all bounded measurable function $\psi$. Since 
$\zeta_i$ commutes with $\eta_j(\psi)= \hat{\eta}_j \psi \hat{\eta}_j$, we verify that 
$\mu(\hat{\eta}_j\zeta_i(\psi)\hat{\eta}_j)=\mu(\hat{\zeta}_i^*\hat{\eta}_j\psi\hat{\eta}_j\hat{\zeta}_i)$ for any measurable partition $(\hat{\eta}_j)$. Hence
$\zeta_i(\psi)=\hat{\zeta}^*_i\psi\hat{\zeta}_i$. We also note that unless we demand that $\hat{\zeta}_i$ are non-negative functions
the family of functions $(\hat{\zeta}_i)$ are not uniquely determined by the family of positive maps $(\zeta_i)$. 
\end{proof} 

\vsp
We denote by $\cll = \{ \ul{\zeta}=(\zeta_i(\psi)=\hat{\zeta}_i^*\psi\hat{\zeta}_i,\;\sum_i\hat{\zeta}_i^*\hat{\zeta}_i=1 \}.$ For any $\ul{\zeta} \in
\cll$ we check that 
$$H_{\mu}(\ul{\zeta}) = -\sum_i \mu(\hat{\zeta}^*_i\hat{\zeta}_i)ln\mu(\hat{\zeta}^*_i\hat{\zeta}_i) + \sum_i\mu(\hat{\zeta}_i^*\hat{\zeta}_iln(\hat{\zeta}^*_i\hat{\zeta}_i))$$
and the following relations as in Kolmogorov-Sinai theory hold by Proposition 2.3:

\vsp
\begin{pro} 
For any finite or countable partitions $\zeta,\eta,\beta \in \cll$ following hold:

\NI (a) $H_{\mu}(\ul{\zeta} \circ \ul{\eta}) \ge H_{\mu}(\ul{\zeta})$

\NI (b) $H_{\mu}(\ul{\zeta} | \ul{\eta} \circ \ul{\beta}) \le H_{\mu}(\ul{\zeta}|\ul{\eta})$

\NI (c) $H_{\mu}(\ul{\zeta} \circ \ul{\eta} | \ul{\beta}) \le H_{\mu}(\ul{\zeta} | \ul{\beta}) + H_{\mu}(\ul{\eta} | \ul{\beta})$

\NI (d) $H_{\mu}(\theta(\ul{\zeta}) | \theta (\ul{\beta})) =H_{\mu}(\ul{\zeta} | \ul{\beta})$
\end{pro} 

\vsp
\begin{proof} 
(a) and (b) follows from the general case. Since $\cll$ is a commutative class, (b) is equivalent
to (c). (d) also follows from the general case. 
\end{proof} 

\vsp
We define dynamical entropy $h_{\mu}(\ul{\zeta},\theta)$ as in Section 3 by equation (16) 
and note in the present case $\ul{\zeta}$ is an invariant partition for $\mu$. 
 
\vsp
\begin{pro} 
For two finite partitions $\zeta,\eta \in \cll$
$$h_{\mu}(\ul{\zeta},\theta) \le h_{\mu}(\ul{\eta},\theta)+H_{\mu}(\ul{\zeta}|\ul{\eta})$$
\end{pro} 

\vsp
\begin{proof} 
In spirit we follow [Pa]. Take $\ul{\zeta}_{(n)}=\theta^{n-1}(\ul{\zeta}) \circ ...\circ \ul{\zeta}$ for 
any partition $\ul{\zeta}$ and verify the following step as in Proposition 3.8
$$H_{\mu}(\ul{\zeta}_{(n)}) \le H_{\mu}(\ul{\eta}_{(n)} \circ \ul{\zeta}_{(n)})$$ 
$$=H_{\mu}(\ul{\eta}_{(n)}) + H_{\mu}(\ul{\zeta}_{(n)}|\ul{\eta}_{(n)})$$ 
$$\le H_{\mu}(\ul{\eta}_{(n)})+ \sum_{0 \le k \le n-1} H_{\mu}(\theta^k(\ul{\zeta})|\theta^k(\ul{\eta}))$$
$$\le H_{\mu}(\ul{\eta}_{(n)}) + nH_{\mu}(\ul{\zeta}|\ul{\eta}).$$
Thus the result follows. 
\end{proof} 

\bigskip
In case $\eta$ is
a simple partition of measurable sets then 
$$H_{\mu}(\ul{\zeta}|\ul{\eta})=
-\sum_{i \in \zeta}\mu(\IE_{\eta}[\tilde{\zeta}_i(1))] ln \IE_{\eta}
[\tilde{\zeta}_i(1)]) + \sum_{i \in \zeta} \mu(\tilde{\zeta}_i(1)ln \tilde{\zeta}_i(1)),$$
where $E_{\eta}(\psi)=\sum_j{<\eta_j, \psi>_{\mu} \over
\mu(\eta_j)} \chi_{\eta_j}$ is the conditional expectation on the $\sigma$-field
generated by the partition $\eta$. So given any partition $\ul{\zeta} \in \cll$
and $\epsilon > 0$ there exists a simple partition $\ul{\eta}$ so that $H_{\mu}(\ul{\zeta}|\ul{\eta})
\le \epsilon.$ Thus a simple consequence of Proposition 4.3, we have the following result.

\vsp
\begin{thm}  
Let $(\Omega,\clf,\mu)$ be a probability space and
$\theta$ is a one-one and onto measurable $\mu$-invariant map. Then
$$\mbox{sup}_{\zeta \in \cll}h_{\mu}(\zeta,\theta)=\mbox{sup}_{\zeta \in \cll_0}
h_{\mu}(\zeta,\theta).$$
\end{thm} 

\vsp
In the following we show that the dynamical entropy is same even when we restrict to the class $\cll \cap C(\Omega)$ (in case $\Omega$ is a locally compact topological space ) for a regular measure $\mu$. 
To that end we start with the following proposition where $\cll_0$ is the class of partition of unity into disjoint measurable sets $(\zeta_i)$. 

\vsp
\begin{pro} 
Let $\Omega$ be a compact Hausdorff space with a regular measure $\mu$. For any fixed 
partition $\zeta \in \cll_0$ and $\epsilon > 0$ there exists a partition $\eta \in \cll \cap C(\Omega)$ so that
$$H_{\mu}(\zeta|\eta) \le \epsilon$$
\end{pro} 

\vsp
\begin{proof} 
Given an element $\zeta \in \cll_0$ i.e. a partition of unity into disjoint measurable 
sets $(\zeta_i)$, we choose by
regularity of the measure $\mu$ a sequence of  
non-negative continuous functions $(\eta^n_i)$ so that $\mu(\eta^n_i\zeta_j) \raro
\mu(\zeta_i)\delta^i_j$. Now we set $\zeta^n_i={\eta^n_i \over \sum_{i \zeta}\eta_i } $. Thus $\zeta^n$ is 
a partition of unity with elements in $C(\Omega)$.  
Thus $0 \le H_{\mu}(\zeta \circ \zeta^n)-H_{\mu}(\zeta^n) 
\le H^c_{\mu}(\zeta \circ 
\zeta^n) - H^c_{\mu}(\zeta^n) \raro 0 $. 
\end{proof} 

\vsp
\begin{thm} 
Let $\Omega$ be a compact Hausdorff space with a 
regular probability measure $\mu$ and $\theta$ be a homeomorphism on 
$\Omega$. Then 
$$h_{\mu}(\theta)=\mbox{sup}_{\zeta \in \cll \cap C(\Omega) } 
h(\zeta,\theta).$$
\end{thm} 

\vsp
\begin{proof} 
It follows from Proposition 6.5 and the basic inequality in Proposition 6.3. 
\end{proof} 

\vsp 
So far we have considered essentially a class of partitions of unity that does not depend on the automorphism $\theta$ under consideration. Our main result in the last section says that $h_{\mu}(\theta)$ remains same as along as there is a measurable partitions $(\ul{\zeta}=(\zeta_i)$ that generates the Borel $\sigma-$field by $\theta$ and we take sup over all possible 
admissible partitions of unity in $\clm$ for $(\clm,\theta,\phi_{\mu})$, where $\clm=L^{\infty}(\Omega,\clf,\mu)$ and $\phi_{\mu}(f)=\int fd\mu$ for $f \in \clm$. 

\vsp
\section{ Translation invariant state in quantum spin chain and its quantum dynamical entropy: }

\vsp 
In this section we include our main application of our main theorem proved in section 3. Let 
$\IM=\otimes_{k \in \IZ} \!M^{(k)}_d(\IC)$ be the two sided quantum spin chain $C^*$-algebra and
$\theta$ be the right translation automorphism on $\IM$ with an invariant state $\omega$ of $\IM$. 
We recall standard definition that two states $\omega_1$ and $\omega_2$ of a $C^*$-algebra, here $\IM$ are called 
{\it quasi-equivalent} if $\pi_{\omega_1}(\IM)''$ and $\pi_{\omega_2}(\IM)''$ are isomorphic i.e. there exists 
an isomorphism $\alpha:\pi_{\omega_1}(\IM)'' \raro \pi_{\omega_2}(\IM)''$ such that $\alpha(\pi_{\omega}(x)) = 
\pi_{\omega_2}(x)$ for all $x \in \IM$. Two translation invariant dynamics $(\IM,\theta,\omega_1)$ and $(\IM,\theta,\omega_2)$ 
are called {\it isomorphic} if there exists an automorphism $\alpha$ on $\IM$ such that $\alpha \theta= \theta \alpha$ and 
$\omega_1 \alpha =\omega_2$. In [Mo5], we proved that Connes-St\o rmer dynamical entropy $h_{CS}(\IM,\theta,\omega)$ is equal to mean entropy $s(\omega)$. This result gave an indirect proof that mean entropy $s(\omega)$ is an invariance i.e. $s(\omega_1)=s(\omega_2)$ for two translation invariant states $\omega_1$ and $\omega_2$ if $\omega_2= \omega_1 \alpha$ for an automorphism $\alpha$ that commutes with $\theta$. For more technical details and additional results, we refer to [Mo5]. This results in [Mo5] used crucially the fact that $\IM$ is a simple $C^*$-algebra. However, this invariance and its isomorphism result did not capture the full generality of Kolmogorov-Sinai-Ornstein theory [Si,Or1] and thus demands an extension to $W^*$-dynamical systems.   

\vsp 
Let $(\clh_{\omega},\pi_{\omega},\zeta_{\omega})$ be GNS space of a translation invariant state $\omega$ of $\IM$ and $\theta:\pi_{\omega}(\IM)'' \raro \pi_{\omega}(\IM)''$ be the induced automorphism extending $\pi_{\omega}(x) \raro \pi_{\omega}(\theta(x))$ for all $x \in \IM$. Thus the normal state $\omega(X)=\langle \zeta_{\omega}, X \zeta_{\omega} \rangle,\; X \in \pi_{\omega}(\IM)''$ is invariant for the automorphism $\theta$ of $\pi_{\omega}(\IM)''$. We define dynamical entropy $h_{\omega}(\theta)$ of the $C^*$-dynamical system $(\IM,\theta,\omega)$ as $h_{\omega}(\theta)$ of the $W^*$-dynamical system $(\pi_{\omega}(\IM)'',\theta,\omega)$.  
Alternatively, we may define $h_{\omega}(\theta)$ to be equal to $h_{\bar{\omega}}(\bar{\theta})$, where $\bar{\theta}$ is the automorphism on the universal enveloping algebra $(\bar{\IM},\bar{\pi})$ defined as at the end of section 3.   

\vsp 
Let $e=(e_i)$ be an orthonormal basis for $\IC^d$ and $D^{(k)}_e(\IC)$ be the algebra of diagonal matrices in $\!M^{(k)}_d(\IC)$ with respect to $(e_i)$. Let $\ID^e = \otimes_{k \in \IZ} D^{(k)}_e(\IC)$ be the maximal abelian $C^*$-sub algebra and $\IE^e_{\omega_0}$ be the norm one projection onto $\ID^e$ preserving the unique normalised trace $\omega_0$ of $\IM$. In the following proposition we recall Proposition 5.4 in [Mo5].  

\vsp 
\begin{pro} 
Let $\omega$ be a translation invariant state of $\IM$. Then there exists an automorphism 
$\alpha$ on $\IM$ commuting with $\theta$ such that 
$$\omega \alpha = \omega \alpha \IE^e_{\omega_0}$$
\end{pro} 

\vsp 
\begin{pro} 
Let $\omega$ be a translation invariant faithful state of $\IM$ such that 
\be 
\omega = \omega \IE^e_{\omega_0}
\ee
Then 
\be 
\pi_{\omega}(\IM)'' \bigcap \pi_{\omega}(\ID^e)' = \pi_{\omega}(\ID^e)''
\ee
\end{pro} 

\vsp 
\begin{proof} 
In case $\omega$ is the normalized trace $\omega_0$, Lemma 3.1 in [Mo5] includes a proof. For the general case, we need a local version of M. Takasaki's theorem [Ta1]
proved in [AC]. Let $\Lambda$ be a finite subset of $\IZ$ and $\IE^{\omega}_{\Lambda}:\pi_{\omega}(\IM) \raro \pi_{\omega}(\IM_{\Lambda})''$ be the generalised $\omega$ conditional expectation. Since the modular group $(\sigma^t_{\omega})$ keeps all elements in $\pi_{\omega}(\ID^e)''$ invariant by (45) and a Theorem of M. Takesaki [Ta1], we get 
$$\IE^{\omega}_{\Lambda}(zxy)=z\IE^{\omega}_{\Lambda}(x)y$$ 
for all $z,y \in \pi_{\omega}(\ID^e_{\Lambda})''$ and $x \in \IM$. In particular, this shows that $\IE^{\omega}_{\Lambda}(x) \in \pi_{\omega}(\ID^e_{\Lambda})'' \bigcap \pi_{\omega}(\ID^e_{\Lambda})' = \pi_{\omega}(\ID^e_{\Lambda})'' $ for all $x \in \pi_{\omega}(\IM)'' \bigcap \pi_{\omega}(\ID^e)'$. Since $\IE^{\omega}_{\Lambda}(x) 
\raro x$ as $\Lambda \uparrow \IZ$, we conclude that $x \in \pi_{\omega}(\ID^e)''$.   
\end{proof} 

\vsp 
Let $\hat{\omega}$ is the coupling state on $\hat{\clm}=\overline{\pi_{\omega}(\IM)'' \otimes \pi_{\omega}(\IM)''}$ defined as in (20) and $s(\hat{\omega},\hat{\omega} \IE^e_{\omega_0} \otimes \IE^e_{\omega_0})$ be the mean relative entropy [OP] defined as in Proposition 3.7 
\be 
s(\hat{\omega},\hat{\omega} \IE^e_{\omega_0} \otimes \IE^e_{\omega_0}) = \mbox{lim}_{n \raro \infty}{1 \over n} S(\hat{\omega}|\hat{\clm}_n,\hat{\omega} \IE^e_{\omega_0} \otimes \IE^e_{\omega_0}|\hat{\clm}_n) 
\ee
where $\hat{\clm}_n= \overline{\pi_{\omega}(\IM^{\IZ_n})'' \otimes \pi_{\omega}(\IM^{\IZ_n})''}$ and $\IM^{\IZ_n} = \otimes_{1 \le  k \le n} \!M^{(k)}_d$. 

\vsp 
\begin{pro} 
Let $\omega$ be a translation invariant faithful state of $\IM$ satisfying (45). Then 
\be 
s(\hat{\omega},\hat{\omega} \IE^e_{\omega_0} \otimes \IE^e_{\omega_0}) \le s(\omega) \le s(\omega_0)=ln(d)
\ee
\end{pro} 

\vsp 
\begin{proof} 
Since $\omega_0 \otimes \omega_0 \IE^e_{\omega_0}= \omega_0 \otimes \omega_0$ on $\hat{\clm}_n$, we have by Theorem 5.15 in [OP]
\be 
S(\hat{\omega}|\hat{\clm}_n,\omega_0 \otimes \omega_0 |\hat{\clm}_n) =
S(\hat{\omega}|\hat{\ID^e}_n,\omega_0 \otimes \omega_0 |\hat{\ID^e}_n) 
+ 
S(\hat{\omega}|\hat{\clm}_n,\hat{\omega} \IE^e_{\omega_0} \otimes \IE^e_{\omega_0}|\hat{\clm}_n) 
\ee
where $\hat{\ID^e}_n= \ID^e_n \otimes \ID^e_n$ and $\ID^e_n= \otimes_{1 \le k \le n} \!D_d^{(k)}(\IC)$. 
Thus 
$$S(\hat{\omega}|\hat{\clm}_n,\hat{\omega} \IE^e_{\omega_0} \otimes \IE^e_{\omega_0}|\hat{\clm}_n)
= S(\hat{\omega}|\hat{\clm}_n) - S(\hat{\omega}|\hat{\ID^e}_n)$$ 
We compute the following identities:
$$S(\hat{\omega}|\hat{\ID^e}_n) $$
$$= S(\omega|\ID^e_n) $$
(since $\hat{\omega}(|e_I><e_I| \otimes |e_J><e_J|)= \delta^I_J \omega(|e_I><e_I|)$ for all $|I|,|J| \le  \infty $ )
$$= S(\omega|\IM^{\IZ_n})$$
(since $\omega \IE^e_{\omega_0} = \omega$ and $\omega_0 \IE^e_{\omega_0} = \omega_0$ as well.
Thus we get 
$$S(\hat{\omega}|\hat{\clm}_n,\hat{\omega} \IE^e_{\omega_0} \otimes \IE^e_{\omega_0}|\hat{\clm}_n)
\le 2S(\omega|\IM^{\IZ_n}) - S(\omega|\IM^{\IZ_n}) = S(\omega|\IM^{\IZ_n})$$
where we have used a well known formula $S(\omega_{12}) \le S(\omega_1) + S(\omega_2)$ for a state $\omega_{12}$ on a bipartite system of $C^*$-algebras with their marginals equals to $\omega_1$ and $\omega_2$ respectively ).
Now using the upper bound in the limit as $n \raro \infty$ in (47), we conclude the proof. 
\end{proof}

\vsp 
\begin{thm} 
Let $\omega$ be an infinite tensor product state of $\IM$. Then $h_{\omega}(\theta)=s(\omega)$.
\end{thm} 

\vsp 
\begin{proof} 
We assume first that $\omega$ is a faithful state on $\IM$. A measurement $\ul{\zeta}=(\zeta_i)$ is called local if $\zeta_i(x)=x$ for all $x \in \IM_{\Lambda'}$, where $\Lambda$ is a finite subset of $\IZ$ and $\Lambda'$ is its complement in $\IZ$. For any $\omega$-invariant local partition $\ul{\zeta}$ of unity, we clearly have $h_{\omega}(\theta,\ul{\zeta}) \le s(\omega)$. For a quick proof, we fix any $\epsilon > 0$ and find $n_{\epsilon} \ge 1$ such that 
$$h_{\omega}(\theta,\ul{\zeta})-\epsilon \le {1 \over n_{\epsilon}}H_{\omega}(\ul{\zeta}_{(n_{\epsilon})}) \le {1 \over n_{\epsilon}} S(\omega|\Lambda_n)$$
where $\Lambda_n= \bigcup_{1 \le k \le n} \Lambda +k$. By our suitable choice for $n_{\epsilon}$, we may assume as $\epsilon \raro 0$, $n_{\epsilon} \raro \infty$. Thus we conclude that $h_{\omega}(\theta,\ul{\zeta}) \le s(\omega)$ by taking $\epsilon \raro 0$. For a general admissible measurement $\ul{\zeta}$, we will approximate by local measurements as follows: For a finite subset $\Lambda$ of $\IZ$, let $\IE_{\Lambda}$ be the norm one projection from $\clm$ onto $\pi_{\omega}(\IM_{\Lambda})''$ preserving the state $\omega$. Since $\omega$ is an infinite tensor product state, such a norm projection exists by a theorem of M. Takesaki [Ta1]. Now we consider the partition of unity defined by 
$\ul{\zeta}_{\Lambda}= (\IE_{\Lambda} \zeta_i \IE_{\Lambda} \otimes I_{\Lambda'})$. That it is also $\omega$-invariant is evident since $\omega$ is an infinite-tensor state. By our argument above, we have $h_{\omega}(\theta,\ul{\zeta}_{\Lambda}) \le s(\omega)$. Now taking sup over all possible finite subset $\Lambda$, we conclude the proof since the map $\ul{\zeta} \raro h_{\omega}(\theta,\ul{\zeta})$ is upper semi-continuous for each fixed $\omega$ as well, where we equip $(\zeta_i)$--Arveson's bounded weak topology [Pau] to ensure $\phi \zeta_i^{\Lambda}(X) \raro \phi \zeta_i(X)$ as $\Lambda \uparrow \IZ$ for all $X \in \pi_{\omega}(\IM)''$. 

\vsp 
Now we drop faithfulness assumption on $\omega$ as in last theorem by taking $\omega_{\lambda}=\lambda \omega + (1-\lambda)\omega_0$ and $h_{\omega_{\lambda}}(\theta,\ul{\zeta}) \le s(\omega_{\lambda})=\lambda s(\omega)+(1-\lambda)s(\omega_0)$ for $\lambda \in (0,1)$ for any admissible partition of unity $\ul{\zeta}$. We use upper-semi-continuity property of the map $\phi_{\lambda} \raro h_{\omega_{\lambda}}(\theta,\ul{\zeta})$ to complete the proof. 
\end{proof}

\vsp 
\begin{thm} 
Let $\omega$ be a translation invariant state of $\IM$. Then  
\be 
s(\omega) \le h_{\omega}(\theta) \le 2s(\omega) 
\ee
\end{thm} 

\vsp 
\begin{proof} 
Let $\omega$ be as well faithful for the time being. Since $h_{\omega}(\theta)$ and $s(\omega)$ are invariances for the dynamics $(\IM,\theta)$, by Proposition 5.1 we may assume without loss of generality that $\omega = \omega \IE^e_{\omega_0}$ holds. We consider the family $\hat{\eta}_j=|e_j \rangle \langle e_j|$ of projections. The partition of unity $\ul{\eta}=(\eta_i)$ is admissible and maximal for $(\pi_{\omega}(\IM)'',\theta,\omega)$ since $\pi_{\omega}(\ID^e)''$ is maximal abelian in 
$\pi_{\omega}(\IM)''$ by Proposition 5.2. 

\vsp 
We consider the $W^*$-dynamics $(\bar{\IM},\bar{\theta},\bar{\omega})$, where $\bar{\omega}_0$ is unique normal extension of the unique tracial state $\omega_0$ of $\IM$ to its universal von-Neumann algebra $\bar{\IM}$. Though $\bar{\omega}_0$ is not faithful on $\bar{\IM}$, a modified argument restricting to its support projection ensures an increasing von-Neumann sub-algebras $\cln_{\hat{\ul{\eta}}_{(n)}}=\vee_{0 \le k \le n-1} \bar{\theta}^{k}(\pi_{\bar{\omega}_0}(\!M^{(k)}_d))$ and faithful normal norm one projections $\IE_{\cln_{\hat{\ul{\eta}}_{(n)}}}$ from $\bar{\pi}(\IM)''$ onto $\cln_{\hat{\ul{\zeta}}_{(n)}}$ preserving $\bar{\omega}_0$ following the construction as in Proposition 3.7. Thus the maximal partition $\ul{\eta}$ satisfies the criteria in Theorem 3.11 with $\bar{\omega}_0$ for the dynamics $(\bar{\IM},\bar{\theta},\bar{\omega})$ and $h_{\bar{\omega}}(\bar{\theta},\ul{\eta})=s(\omega)$. Thus by Theorem 3.11 and Proposition 5.3, we complete a proof for (49) for faithful $\omega$.  

\vsp 
Now we will drop our assumption on $\omega$. For a general state $\omega$, we consider the faithful state $\omega_{\lambda}=\lambda \omega + (1-\lambda)\omega_0$ for $0 \le \lambda < 1$. Then by the first part of our argument we have for any admissible partition $\ul{\zeta}$ for $(\clm,\theta,\omega_{\lambda})$:
$$h_{\omega_{\lambda}}(\theta,\ul{\zeta}) \le 2s(\omega_{\lambda}),$$ 
where $s(\omega_{\lambda}) = \lambda s(\omega) + (1-\lambda)s(\omega_0)$. Since $\omega_{\lambda} \raro h_{\omega_{\lambda}}(\theta,\ul{\zeta})$ is an 
upper-semi-continuous function for each fixed $\ul{\zeta}$, we get $h_{\omega}(\theta,\ul{\zeta}) \le 2s(\omega)$ by taking limit $\lambda \raro 1$. That $h_{\omega}(\theta) \ge s(\omega)$ follows as $h_{\omega}(\theta) \ge h_{\omega}(\theta,\ul{\eta})=s(\omega)$. This completes the proof. 
\end{proof} 

\bigskip
{\centerline {\bf REFERENCES}}

\begin{itemize} 
\bigskip 
\item{[Ar]} Araki, H.: Relative entropy of states of von-Neumann algebra,
Publ. RIMS, Kyoto Univ., 11, pp. 809-833, 1976.

\item{[CS]} Connes, A.; Størmer, E.: Entropy of automorphisms of II$_1$ -von Neumann algebras, Acta Math. 134 (1975), 289-306.

\item{[CNT]} Connes, A., Narnhofer, H. and Thirring, W.: Dynamical entropy of $C^*$-
algebras and von Neumann algebras, Commun. Math. Phys. 112 (1987), 691-719.

\item{[CFS]} Cornfeld, I. P. Fomin, S. V.; Sinaĭ, Ya. G. Ergodic theory. Translated from the Russian by A. B. Sosinskiĭ. Grundlehren der Mathematischen Wissenschaften [Fundamental Principles of Mathematical Sciences], 245. Springer-Verlag, New York, 1982.

\item {[Ka]} Kadison, Richard V.: A generalized Schwarz inequality and algebraic invariants for operator algebras,  Ann. of Math. (2)  56, 494-503 (1952). 

\item {[Mo1]} Mohari, A.: Pure inductive limit state and Kolmogorov property. II 
Journal of Operator Theory. vol 72, issue 2, 387-404 (2014).   
    
\item {[Mo2]} Mohari, A.: Translation invariant pure state on $\otimes_{k \in \IZ}\!M^{(k)}_d(\IC)$ and Haag duality, Complex Anal. Oper. Theory 8 (2014), no. 3, 745-789.

\item{[Mo3]} Mohari, A.: Translation invariant pure state on $\clb=\otimes_{k \in \IZ}\!M^{(k)}_d(\IC)$ and its split property, J. Math. Phys. 56, 061701 (2015).

\item{[Mo4]} Mohari, A.: Isomorphism theorem for Kolmogorov states of $\IM=\dsp{\otimes_{n \in \IZ}}\!M^{(n)}_d(\IC)$, arXiv:1309.7606

\item{[Mo5]} Translation invariant state and its mean entropy-I arXiv:1705.11038. 
    
\item{[Mo6]} Mohari, A.: Translation invariant states in free probability; in preparation. 

\item {[NS]} Neshveyev, S.: St\o rmer, E. : Dynamical entropy in operator algebras. Springer-Verlag, Berlin, 2006. 

\item {[Ne]}  Neumann, J. von: Mathematische Grundlagen der Quantenmechanik, Springer, Berlin, 1932.

\item {[Or1]} Ornstein, D. S.: Two Bernoulli shifts with infinite entropy are isomorphic. Advances in Math. 5 1970 339-348 (1970).  

\item{[Or2]} Ornstein, D. S.: A K-automorphism with no square root and Pinsker's conjecture, 
Advances in Math. 10, 89-102. (1973).

\item {[OP]} Ohya, M., Petz, D.: Quantum entropy and its use, Text and monograph in physics, Springer-Verlag 1995. 

\item {[Pau]} Paulsen, V.: Completely bounded maps and operator algebras, Cambridge Studies in Advance Mathematics 78, Cambridge University Press. 2002

\item {[Pa]} Parry, W.: Topics in Ergodic Theory, Cambridge University Press, 1981. 

\item {[Sak]} Sakai, S.: C$^*$-algebras and W$^*$-algebras, Springer 1971.  

\item {[Si]} Sinai, Ja.: On the concept of entropy for a dynamic system. (Russian) Dokl. Akad. Nauk SSSR 124 1959 768-771. 
 
\item {[SS]} Sinclair, Allan M.; Smith, Roger R.: Finite von Neumann algebras and masas. London Mathematical Society Lecture Note Series, 351. Cambridge University Press, Cambridge, 2008.

\item{[St\o 2]} St\o rmer, E.: A survey of non-commutative dynamical entropy. Classification of nuclear $C^*$-algebras. Entropy in operator algebras, 147-198, Encyclopaedia Math. Sci., 126, Springer, Berlin, 2002.

\item{[Ta1]} Takesaki, M.: Conditional Expectations in von Neumann Algebras, J. Funct. Anal., 9, pp. 306-321 (1972)
 
\item{[Ta2]} Takesaki, M. : Theory of Operator algebras II, Springer, 2001.
  
\end{itemize}

\end{document}